\RequirePackage{fix-cm}
\documentclass[smallextended]{svjour3}       
\smartqed  
\usepackage{amsmath,amsfonts}
\usepackage[pagewise, modulo]{lineno}
\usepackage[normalem]{ulem}
\usepackage{graphicx}
\usepackage[percent]{overpic}
\usepackage[caption=false]{subfig}
\usepackage{epsfig} 
\usepackage{epstopdf}
\usepackage[colorlinks=true]{hyperref}

\usepackage[ruled,vlined]{algorithm2e}

\newcommand{\trp}{^{\text{\scriptsize T}} }

\newcommand{\inv}{^{-1} }
\newcommand{\sbr}[1]{\ensuremath{_{\mathrm{#1}}}}
\newcommand{\spr}[1]{\ensuremath{^{\mathrm{#1}}}}

\def\diag{\operatorname{diag}}
\begin{document}


\title{Randomized maximum likelihood based posterior sampling
}

\titlerunning{RML based posterior sampling}        

\author{Yuming Ba        \and
        Jana de Wiljes   \and
        Dean S. Oliver \and
        Sebastian Reich
}

\authorrunning{Y. Ba, J. de Wiljes, D. S. Oliver and S. Reich} 

\institute{Yuming Ba \at
              School of Mathematics, Hunan University, Changsha 410082, China \\
              \email{yumingb9630@163.com}             \\
             \emph{Present address:} School of Mathematics and Systems Science, Guangdong Polytechnic Normal University, Guangzhou 510665 
           \and
          Jana de Wiljes \at
              Universit\"at Potsdam,
Institut f\"ur Mathematik, Karl-Liebknecht-Str.~24/25, D-14476 Potsdam, Germany\\
\email{wiljes@uni-potsdam.de}
         \and
           Dean S. Oliver \at
              NORCE Norwegian Research Centre, Bergen, Norway\\
\email{dean.oliver@norceresearch.no}
\and
          Sebastian Reich \at
              Universit\"at Potsdam,
Institut f\"ur Mathematik, Karl-Liebknecht-Str.~24/25, D-14476 Potsdam, Germany\\
\email{sebastian.reich@uni-potsdam.de}
}

\date{Received: date / Accepted: date}

\maketitle

\begin{abstract}
Minimization of a stochastic cost function is commonly used for approximate sampling in high-dimensional Bayesian inverse problems with Gaussian prior distributions and multimodal posterior distributions. The density of the samples generated by minimization is not the desired target density, unless the observation operator is linear, but the distribution of samples is useful as a proposal density for importance sampling or for Markov chain Monte Carlo methods. In this paper, we focus on applications to sampling from multimodal posterior distributions in high dimensions. We first show that sampling from multimodal distributions is improved by computing all critical points instead of only minimizers of the objective function. For applications to high-dimensional geoscience inverse problems, we demonstrate an efficient approximate weighting that uses a low-rank Gauss-Newton approximation of the determinant of the Jacobian. The method is applied to two toy problems with known posterior distributions and a Darcy flow problem with multiple modes in the posterior.

\end{abstract}


\section{Introduction}

In several fields, including groundwater management, groundwater remediation, and petroleum reservoir management, there is a need to characterize permeable rock bodies whose properties are spatially variable. In most cases, the reservoirs are deeply buried, the number of parameters needed to characterize the porous medium is large and the observations are sparse and indirect \cite{oliver:11}.
In these applications, the problem of estimating model parameters is almost always underdetermined and the desired solution is not simply a best estimate, but rather a  probability density on model parameters conditioned to the observations and to the prior knowledge \cite{tarantola:05}. Because the model dimension in geoscience applications is always large, the posterior distribution is often represented empirically by samples from the posterior distribution.

Unfortunately, the posterior probability density for reservoir properties, conditional to rate and pressure observations, is typically complex and not easily sampled. Several authors have shown that the log-posterior for subsurface flow problems is not convex in some situations \cite{zhang:03b,tavassoli:05,oliver:18a} so that Gaussian approximations of the posterior distribution appear to be dangerous.
Markov chain Monte Carlo methods (MCMC) are often considered to provide the gold standard for sampling from the posterior. It is frequently suggested to be the method against which other methods are compared \cite{law:12}, yet it can be difficult to design efficient transition kernels \cite{efendiev:06,cui:16} and convergence of MCMC can be difficult to assess \cite{gelman:96b}. The number of likelihood function evaluations required to obtain a modest number of independent samples may be excessive for highly nonlinear flow problems \cite{oliver:96d}.

Importance sampling methods can also be considered to be exact sampling methods as they implement Bayes rule directly. They are very difficult to apply in high dimensions, however, as an efficient implementation requires a proposal density that is a good approximation of the posterior \cite{bengtsson:08,mackay:03}. Various methods have been developed in the data assimilation community to ensure that particles are located in regions of high probability density \cite{vanleeuwen:19}.
Although not introduced as importance sampling approaches,
a variety of methods based on minimization of a stochastic objective function have been developed, beginning with Kitanidis \cite{kitanidis:95} and Oliver et al.\ \cite{oliver:96e} who introduced minimization of a stochastic objective function as a way of simulating samples from an approximation of the posterior when the prior distribution is Gaussian and the errors in observations are Gaussian and additive. The distribution of samples based on minimizers of the objective function was shown to be correct for Gauss-linear problems, but when the observation operator $g(m)$ is nonlinear, it was necessary to weight the samples because the sampling was only approximate in that case.

The randomized maximum likelihood (RML) method approach to sampling has been used without weighting in high dimensional inverse problems with Gaussian priors \cite{gao:06a,eydinov:09,cardiff:12,bjarkason:19}. Weights are seldom computed for several reasons: computation of exact weights is infeasible in large dimensions because the   computation of weights requires the second derivative of the observation operator, the proposal density     does not always cover the target density, and sampling without weighting sometimes provides a good approximation of the posterior even in posterior distributions with many modes \cite{oliver:17}. In practice, the most popular implementations are the ensemble-Kalman based forms of the RML method \cite{chen:12a,sakov:12,white:18a,jardak:18} in which a single average sensitivity is used for minimization so that weighting of samples is not possible.

Bardsley et al.\ \cite{bardsley:14} proposed another minimization-based  sampling methodology, randomize-then optimize (RTO), in which the need for computation of the second derivative of the observation operator is avoided for weighting.  The RTO method has then been modified to allow application in very high dimensions \cite{bardsley:20}, but the method is restricted to posterior distributions with a single mode.
Wang et al.\ \cite{wang:18} discussed the relationship between the cost function in the RML method and the cost function in the RTO method, and showed that the methods are equivalent for linear observation operators.
Wang et al.\ also showed that some approximations to the weights in RML could be computed in high dimensions and provided a useful comparison of sampling distributions.

For nonconvex log-posteriors, there could be many (local) minimizers of the RML and RTO cost functions. Oliver et al.\ \cite{oliver:96e} and Oliver \cite{oliver:17} suggested that only the global minimizer should be used for sampling, although finding the global minimizer would be difficult to ensure. To improve the likelihood of converging to the global minimizer, they suggested using  the unconditional sample from the prior as the starting point for minimization. In contrast, Wang et al.\  \cite{wang:18} investigated the effect of various strategies for choosing the initial guess on sample distribution and found that a random initial guess worked well.

Unlike previous methods that compute minimizers only, we show that exact sampling is possible when all critical points of a stochastic objective function are computed and properly weighted.   Computing the weights accurately, however, for all critical points in high dimensions does not appear to be feasible, but we demonstrate that Gauss-Newton approximations of the weights  provide good approximations for minimizers of the objective function in problems with multimodal posteriors. The Gauss-Newton approximations of weights can be obtained as by-products of Gauss-Newton minimization of the objective function, or as low-rank approximations using stochastic sampling approaches.    We also show that valid sampling can be performed without computing all critical points, but by instead randomly sampling of the critical points.

We investigated the performance of both exact sampling and approximate sampling on two small toy problems for which the sampled distribution can be compared with the exact posterior probability density. For a problem with two modes in the posterior pdf, the distribution of samples from weighted RML using all critical points appears to be correct. Approximate sampling using  Gauss-Newton approximation of weights and minimizers samples well from both modes, but under-samples the region between modes. The data misfit is not a useful approximation of the weights in this case.

We also applied the approximate sampling method to the problem of estimating permeability in a 2D porous medium from 25 measurements of pressure. In this case, the distribution of weights was relatively large, even when the log-permeability was distributed as multivariate normal. High-dimensional state spaces such as considered here have a severe effect on importance sampling and remedies such as tempering are suggested to reduce the impact of the dimensionality on the estimation \cite{Beskos:14}.

\section{RML sampling algorithm}

Given a prior Gaussian distribution ${\rm N}(\bar m,C_M)$
on a set of model parameters $m \in {\mathbb R}^{N_m}$
and observations $d^{\rm o} \in \mathbb{R}^{N_d}$ which are related
to the model parameters through a forward map
$g:\mathbb{R}^{N_m}\to \mathbb{R}^{N_d}$ for unknown
$m^\ast$ and unknown measurement errors $\epsilon \sim {\rm N}(0,C_D)$,
i.e.,
\begin{equation}
d^{\rm o} = g(m^\ast) + \epsilon,
\label{eq:obmodel}
\end{equation}
we wish to generate samples $m_i$, $i=1,\ldots,N_e$, from the posterior distribution
\begin{equation}
\pi_M(m|d^{\rm o}) = \frac{\pi_{MD}(m,d^{\rm o})}{\pi_{D}(d^{\rm o})}
\propto \exp (-L(m)) \label{eq:posterior}
\end{equation}
with negative log likelihood function
\begin{equation}
L(m) = \frac{1}{2} \left(m-\bar m\right)^T C_M^{-1}\left(m-\bar m\right) + \frac{1}{2}
\left(g(m)-d^{\rm o} \right)^T C_D^{-1} \left( g(m)-d^{\rm o} \right).
\label{eq:neglog}
\end{equation}
The normalisation constant $\pi_D(d^{\rm o})$ is unknown, in general. We will use $\pi_M(m) := \pi_M(m|d^{\rm o})$ in order to simplify notation.

In this paper, we will show how to use the RML method
in order to produce independent weighted
Monte Carlo samples from the posterior distribution \eqref{eq:posterior}.
Posterior sampling problems of the form \eqref{eq:posterior} with negative log likelihood function
\eqref{eq:neglog} arise from many practical Bayesian inference problems. In practical applications, where the number of model parameters is typically large, the computation of exact weights is infeasible. For those cases we suggest approximations.

%
%

\subsection{The trial distribution: RML as proposal step}

The RML method draws samples $(m_i',\delta_i')$, $i=1,\ldots,N_s$, from the Gaussian
distribution
\begin{multline}
q_{M'\Delta'} (m',\delta') = \frac{1}{(2\pi)^{N_m N_d/2} |C_M|^{1/2} |C_D|^{1/2}} \\
 \times\exp \left( - \frac{1}{2} \left(m' -\bar m\right)^T C_M^{-1}\left(m' -\bar m\right)  - \frac{1}{2}
 \left(\delta' -d^{\rm o}\right)^T C_D^{-1}\left(\delta' -d^{\rm o}\right)  \right)
 \label{eq:xdprop}
\end{multline}
for given $\bar m$ and $d^{\rm o}$ and then computes critical points of the cost functional
\begin{equation}
L_i(m) = \frac{1}{2} \left(m- m_i' \right)^T C_M^{-1}\left(m- m_i' \right) + \frac{1}{2}
\left(g(m)-\delta_i' \right)^T C_D^{-1} \left( g(m)-\delta_i' \right).
\label{eq:Jix}
\end{equation}
by solving
\begin{equation} \label{eq:critical}
 \nabla_m L_i(m) = 0,
\end{equation}
for $m$.
Dropping the subscript $i$, this leads  to a map from $(m,\delta)$ to $(m',\delta')$ defined by
\begin{equation}
\left\{
\begin{aligned}
 &m'  =  m +  C_M G\trp C_D\inv (g(m) - \delta)
  \label{eq:muc}
 \\
 &\delta'    =  \delta
 \end{aligned}
 \right.
 \end{equation}  
which we denote compactly as
\begin{equation}
 z'  =  \Psi(z),
  \label{eq:Psi}
 \end{equation}
 where $z = (m,\delta)$,  $z' = (m',\delta')$  and the differential of $g$ is denoted $G = D g(m)$. 
 The mapping \eqref{eq:Psi} is, in general, not invertible and, hence, a single draw $(m'_i,\delta'_i)$ from (\ref{eq:xdprop}) can lead to multiple critical points $(m_j,\delta_j)$.\footnote{See Sec.~\ref{sec:quadratic} for a sampling problem with a quadratic observation operator, $g$, resulting in a non-invertible mapping. In that example, \eqref{eq:muc} is cubic in the variable $m$. Non-invertible mappings appear to be common for Darcy flow problems (e.g., Sec.~\ref{sec:non-monotonic}).}
We  therefore introduce the set-valued
$$
\mathcal{M}_{z'} = \Psi^{-1} (z')
$$
and denote its elements by $z_j(z') \in \mathcal{M}_{z'}$, $j = 1,\ldots,n(z')$, where $n(z')$ denotes the cardinality of
$\mathcal{M}_{z'}$. 
Each $z$ leads to a unique $z'$, hence the sets $\mathcal{M}_{z'}$ are
disjoint.  Let us
denote the set of all $z'$ for which $n(z')>0$ by $\mathcal{U}$ and let us assume for now that $\mathcal{U}$ agrees
with the support of the distribution $q_{M'\Delta'}$.\footnote{If there are points $z'$ for which $\mathcal{M}_{z'}$ is the empty set, that is, $n(z') = 0$, we  adjust the PDF $q_{M' \Delta'}$ such that $q_{M'\Delta'}(z') = 0$ for $\mathcal{M}_{z'} = \emptyset$.}

A distribution $q(z')$ transforms under a map (\ref{eq:Psi}) into a distribution $p(z)$ according to
\begin{equation}
q(z')  = \sum_{z_i \in \mathcal{M}_{z'}} \frac{p(z_i)}{J(z_i)}
  \label{eq:RML_proposal}
\end{equation}
with Jacobian $J(z) = \det(D\Psi(z))$. We will frequently use the abbreviation $|A|$ for the determinant $\det (A)$ of a matrix $A$. An explicit expression for $p(z)$ is obtained via
\begin{equation*}
p(z) = n(z')^{-1} J(z)\, q(z')
\end{equation*}
for all $z \in \mathcal{M}_{z'}$, which satisfies \eqref{eq:RML_proposal}.
In the original notation and employing (\ref{eq:muc}),
the transformed distribution $p_{M\Delta}$ is given by
\begin{equation*}
\begin{aligned}
p_{M\Delta}(m,\delta) &:= n(m')^{-1}q_{M'\Delta'}(m',\delta') \, J(m,\delta)  \\
&= n(m')^{-1}q_{M'}\left(m +  C_M G\trp C_D\inv (g(m) - \delta) \right) \, q_{\Delta'} \left(  \delta
  \right) \, J(m,\delta)  \, .
\end{aligned}
\end{equation*}
Here $n(m')$ is the total number of critical points of  \eqref{eq:Jix} for each $(m',\delta')$ and $J(m,\delta)$ denotes the Jacobian determinant associated with the  map
$(m,\delta) \to (m',\delta')$. In the following, we assume that $J \not=0$ everywhere, i.e., the map is locally invertible.  The Jacobian matrix is provided by
\begin{equation}
\left( \begin{array}{cc} I +  D b(m,\delta) & -  C_M G^{\rm T} C_D^{-1} \\
0 & I \end{array} \right)
\label{eq:jacobian}
\end{equation}
with $b(m,\delta) = C_M G^{\rm T}C_D^{-1}(g(m)-\delta)$.

Hence, given samples, ($m'_i,\delta'_i)$,  from \eqref{eq:xdprop} we can easily produce samples,
$(m_k,\delta_k)$,
from the distribution $p_{M\Delta}(m,\delta)$ and would like to use them
as importance samples from the target distribution
$\pi_M(m) := \pi_M(m|d^{\rm o})$ as defined by \eqref{eq:posterior}.
Note that the target density $\pi_M(m)$ does not specify a distribution in $\delta$ and we will
explore this freedom in the subsequent discussion in order to define an efficient importance
sampling procedure.

Indeed, we may introduce an extended target distribution by
\begin{equation*}
\pi_{M\Delta}(m,\delta) := \pi_M(m) \, \pi_\Delta(\delta |m)
\end{equation*}
without changing the marginal distribution in $m$.
The conditional distribution $\pi_\Delta(\delta |m)$
will be chosen to make the proposal density similar to the target density, i.e.,
\begin{equation*}
\pi_{M\Delta}(m,\delta) \approx p_{M\Delta}(m,\delta).
\end{equation*}
We will find that equality can be achieved for linear forward maps, $g(m) = Gm$.
In all other cases, samples, $(m_k,\delta_k)$ from $p_{M\Delta}(m,\delta)$
will receive an importance weight
\begin{equation}
w_k \propto \frac{\pi_{M\Delta}(m_k,\delta_k)}{p_{M\Delta}(m_k,\delta_k)}\label{eq:weights}
\end{equation}
subject to the constraint $\sum_{k=1}^{N_e} w_k = 1$. Note that all involved distributions need
only to be available up to normalisation constants which do not depend on $m$ or $\delta$.

The subsequent discussion will reveal a natural choice for $\pi_\Delta(\delta |m)$ and
will lead to an explicit expression for \eqref{eq:weights}. Let us therefore
go through the  analysis to factor
\begin{equation*}
p_{M\Delta}(m,\delta) = n(m')^{-1} q_{M'\Delta'}(m',\delta') \, J(m,\delta)
\end{equation*}
to determine a candidate $\pi_\Delta(\delta |m)$.
First, we expand the negative log density $- \log q_{M'\Delta'}$, ignoring the normalization constant:
\begin{multline*}
  \frac{1}{2} \left(m -\bar m +  C_M G\trp C_D\inv (g(m) - \delta)
  \right)\trp C_M\inv  \left(m -\bar m +  C_M G\trp C_D\inv (g(m) - \delta) \right)\\
 + \frac{1}{2}    \left( \delta
  -d^o \right)\trp C_D\inv \left(  \delta
  -d^o \right) \\
 =   \frac{1}{2} \left(m -\bar m \right)\trp  C_M\inv \left(m -\bar m \right)
   +  \frac{1}{2} (g(m) - \delta)\trp C_d\inv G C_M G\trp  C_D\inv (g(m) - \delta)   \\
 + \frac{1}{2 }  (g(m) - \delta)\trp C_d\inv  G \left(m -\bar m\right)
 + \frac{1}{2 } \left(m -\bar m\right)\trp G\trp C_D\inv (g(m) - \delta)  \\
 + \frac{1}{2}    \left(  g(m) -d^o \right)\trp C_D\inv \left( g(m) -d^o \right)
 + \frac{1}{2 }    \left(  g(m) -\delta \right)\trp C_D\inv \left( g(m) -\delta \right) \\
 - \frac{1}{2 }    \left(  g(m) -\delta \right)\trp C_D\inv \left( g(m) -d^o \right)
 - \frac{1}{2 }    \left(  g(m) -d^o \right)\trp C_D\inv \left( g(m) -\delta \right).
\end{multline*}
To simplify the notation we will use
\begin{equation}  V :=   C_D + G C_M G\trp \label{eq:V} \end{equation}
and
\begin{equation} 
\eta(m) 
 :=   G ( m - \bar m) -( g(m) -d^o)    .
\label{eq:eta}
\end{equation}
Then, using  the new definitions  to simplify notation, we obtain
\begin{equation}  \begin{split}
 &p_{M\Delta}(m,\delta)  = \\
 &\overbrace{A_0\, \exp \left[-\frac{1}{2} \left(m -\bar m \right)\trp  C_M\inv \left(m -\bar m \right)
- \frac{1}{2}    \left(  g(m) -d^o \right)\trp C_D\inv \left( g(m) -d^o \right) \right]}^{\pi_{M}(m) } \\
& \quad \times \overbrace{A_1\, |V|^{1/2} \exp \left[-\frac{1}{2 } \left( \delta- g(m)   -  V\inv \eta(m)  \right)\trp V \left( \delta - g(m) -   V\inv \eta(m)  \right) \right]}^{\pi_{\Delta}(\delta |m) } \\
& \qquad \qquad \times n(m')\inv A_2\, |V|^{-1/2} \exp \left[ \frac{1}{2} \eta(m)\trp V\inv \eta(m) \right] \, J(m,\delta) ,
\end{split}
\label{eq:simplify_qMD}
\end{equation}
where $A_0$, $A_1$, and $A_2$ are all normalisation constants, independent of $m$ and
$\delta$. $A_0$ is determined from the requirement that
$\int \pi_{M}(m)\, {\rm d}m = 1$. Similarly, $A_1$ is determined from the requirement that
$\int \pi_{\Delta}(\delta |m)\, {\rm d}\delta = 1$. Finally, $A_2$ is determined from the requirement that
$\int p_{M\Delta}(m,\delta)\, {\rm d}m \,{\rm d}\delta = 1$. The last line of   \eqref{eq:simplify_qMD}
is exactly the difference between the proposal density and the target density, which determines the
importance weights  \eqref{eq:weights}.

Note that if the observation operator is linear, then $n(m')=1$ and all terms on the last line of  \eqref{eq:simplify_qMD} are independent of $m$ so the target and proposal densities are equal:  $p_{M\Delta} = q_{M\Delta}$.

\subsection{Weighting of RML samples}\label{subsec:WeRML}
To weight the RML samples, we compute the weights by
\[
w_k \propto \frac{\pi_M(m_k) \, \pi_{\Delta}(\delta_i |m_k)}{p_{M\Delta}(m_k,\delta_k)}
\]
with $\pi_\Delta (\delta |m)$ as defined in \eqref{eq:simplify_qMD}. So the weight on a sample is
\begin{equation}
w   \propto n(m')\, |V|^{1/2} \exp \left[ -\frac{1}{2} \eta(m)\trp V\inv \eta(m) \right] \, J\inv (m,\delta) .
\label{eq:weight-RML}
\end{equation}
The Jacobian determinant and the gradient of the misfit term with respect to the parameter are necessary for the computation of weights. For the low-dimensional space, it is easy to calculate them. However, the computation of the Jacobian determinant and the gradient is difficult when the problems are strongly nonlinear.

In this section, we use the low-rank approximation to get the Jacobian determinant and $\det V$. Using the Gauss-Newton approximation for the Jacobian matrix given by \eqref{eq:jacobian}, we have
\[
J(m,\delta)\approx |I+C_MG\trp C_D^{-1}G|
\]
and $J$ becomes independent of $\delta$.
Let $m_\text{MAP}$ and $H_\text{misfit}$ denote the minimizer point of (\ref{eq:critical}) and the Hessian matrix of $L_i(m)$ with respect to the misfit term at $m_\text{MAP}$, respectively.
Thus the Hessian matrix of $L_i(m)$ at $m_\text{MAP}$ is given by
\[
H_\text{map}= C_M^{-1}+H_\text{misfit}.
\]
To compute its determinant, we would like to approximate $H_\text{map}$ with a relatively small number of terms. Thus we solve the following generalized eigenvalue problem (GEP):  find $U\in\mathbb R^{N_m\times N_m}$ and
$\Lambda = \diag(\lambda_i) \in \mathbb{R}^{N_m\times N_m}$, which are the generalized eigenvectors and eigenvalues of the matrix pair $H_\text{misfit}$ and  $C_M\inv$, respectively:
\[
 H_\text{misfit} U = C_M\inv U \Lambda,
\]
such that
\[
U\trp C_M\inv U = I \qquad \text{and} \qquad  H_\text{misfit} = C_M\inv U {\Lambda} U\trp C_M\inv.
\]

For the large-scale flow problem, we consider the  Whittle-Mat{\'e}rn prior covariance operator based on the inverse of an elliptic differential operator,
\begin{equation}
\begin{split}
\mathcal{C}_\text{prior} 
& = (-\gamma\Delta+\alpha I)^{-2}, \\
& = \frac{1}{4 \pi \gamma \alpha} \left( 
\frac{r}{ \sqrt{\gamma/\alpha}} \right) \, K_1 \left( 
\frac{r}{ \sqrt{\gamma/\alpha}} \right)
\end{split}
\label{eq:prior}
\end{equation}
where $K_1$ denotes the modified Bessel function of the second kind of order 1. Eq.~\eqref{eq:prior} provides a sparse representation of the inverse covariance and a square root factorization that is useful for computing a low-rank approximation of the Hessian matrix \cite{bui:13}. From \eqref{eq:prior} the variance is seen to be $(4 \pi \alpha \gamma)\inv$ and the range of the covariance to be proportional to $\sqrt{\gamma/\alpha}$. 
%
$C_M$ is given by the discretization  of $\mathcal{C}_\text{prior}$ and the inverse of $C_M$ can be easily factored
\[
C_M\inv=Q Q\trp .
\]
So we have
\[
H_\text{map}=Q Q\trp U \Lambda U^T Q Q\trp +Q Q\trp =Q\hat U(\Lambda +I){\hat U}\trp Q\trp ,
\]
where $\hat U=Q\trp U$ is the matrix of orthonormal eigenvectors for $Q\inv H_\text{misfit}Q^{-{\rm T}}$. We actually want the determinant
\[
|C_M H_\text{map}|=|Q\inv H_\text{map} Q^{-\trp}|=\Big|\hat U(\Lambda +I){\hat U}\trp\Big|=\prod_{i=1}^{N_m}(1+\lambda_i).
\]
When the generalized eigenvalues $\{\lambda_i\}$ decay rapidly, we can use a low-rank approximation of $H_\text{misfit}$ by retaining only the $r$ largest eigenvalues and corresponding eigenvectors, i.e.,
\[
H_\text{misfit}\approx Q Q\trp U_r \Lambda_r U_r\trp Q Q\trp .
\]
Thus we have
\[
J(m,\delta)\approx\prod_{i=1}^r(1+\lambda_i).
\]

We also need the determinant of $V$ for the computation of weights in \eqref{eq:weight-RML}. In section \ref{sec:flow}, we will use a diagonal matrix for $C_D$, i.e., $C_D=\sigma^2 I$. Due to $H_\text{misfit}\approx G\trp C_D\inv G$, we have
\[
G\trp  \sigma^{-2}I G\approx C_M\inv U {\Lambda} U\trp C_M\inv.
\]
Then
\[
C_M G\trp G\approx \sigma^2 U {\Lambda} U\trp C_M\inv.
\]
The determinant of $V$
\[
\begin{aligned}
|V|=|\sigma^2 I+G C_M G\trp |&=\sigma^{2(N_d-N_m)}|\sigma^2 I+C_M G\trp G|\\
&\approx \sigma^{2(N_d-N_m)}|\sigma^2 I+\sigma^2 U {\Lambda} U\trp C_M\inv|\\
&=\sigma^{2 N_d}|I+Q\trp U {\Lambda} U\trp Q|\\
&=\sigma^{2 N_d}|I+\hat U {\Lambda} \hat U\trp|\\ &\approx \sigma^{2 N_d}\prod_{i=1}^r(1+\lambda_i).\\
\end{aligned}
\]
Thus the determinant of $V$ can be also replaced by a low-rank approximation. $\sigma^{2N_d}$ is not necessary in $|V|$ because it appears in all weights and can be factored out. Then \eqref{eq:weight-RML} can be approximated by
\[
w\propto n(m')\, \exp \left[ -\frac{1}{2} \eta(m)\trp V\inv \eta(m) \right] \,\prod_{i=1}^r(1+\lambda_i)^{-1/2}.
\]

In the approach described above, the computation of $G$ for the weights is necessary.
To obtain $G$ for the flow problem in section \ref{sec:flow}, we solve an adjoint system \cite{villa:21}.
In section \ref{sec:quadratic}, we investigate the effect of a Gauss-Newton approximation for the weights for cases in which all the critical points and only minimizers of (\ref{eq:critical}) are obtained, respectively.  As the Gauss-Newton approximation of the weights is shown to be poor for  maximizers of the objective function, we only consider the minimizers in the Darcy flow example (Sec.~\ref{sec:flow}).


\subsection{Weighted RML sampling algorithm}

In this section, we consider two possible situations when seeking independent samples from \eqref{eq:posterior} for a log-posterior  of the form \eqref{eq:neglog}. In both cases, we allow for the possibility that the stochastic cost function  \eqref{eq:Jix}  is nonconvex. Note that in this case the number of critical points may be greater than 1. In the first algorithm, we assume that all critical points can be identified, while in the second case, we suppose that it is only possible to identify a single critical point, but that the total number of critical points is unknown.

\subsubsection{All critical points found}

For problems in low dimensions, with polynomial observation operators $g$, or convex cost functions,  it may be feasible to find all critical points for each pair $(m',\delta')$. If the $i$th sample of $(m',\delta')$, generates a cost function $L_i(m)$ with $n_{c_i}$  critical points, it is possible to sample correctly by weighting each critical point  using \eqref{eq:weight-RML}.

\vspace{\baselineskip}
\begin{algorithm}[H]
 \SetAlgoLined
 \caption{Weighted RML -- computing all the critical points}
	\begin{itemize}
	\item[] $i=1$, $k=1$
	\item[] while $k\le N_e$
\begin{itemize}
		\item[] generate samples $m_i'$	 and $\delta_i'$ from $q_{M'\Delta'} (m',\delta') $
		\item[]  for $j = 1$ to  $n_{c_i}$
		\begin{itemize}
			\item[] solve (\ref{eq:critical}) for $m^{(j)}$ and set $\delta^{(j)}=\delta_i'$
			\item[] compute $w^{(j)}$ using \eqref{eq:weight-RML}
			\item[] assign $m_k = m^{(j)}$ and $w_k' = w^{(j)}$
			\item[] $k=k+1$
		\end{itemize}
		\item[] $i=i+1$
	\end{itemize}
\item[] assign $w_{k} = w'_{k} /\sum_k w'_{k}$
\end{itemize}\label{Alg1}
\end{algorithm}
\vspace{\baselineskip}

Note that when the forward operator $g$ is linear, i.e., $g(m)=Gm$, the stochastic cost function is convex for each pair $(m',\delta')$. The log-posterior has a single critical point, which is the minimizer. Thus, for linear forward operators, we just need to compute the minimizers and the weights are all equal by using  \eqref{eq:weight-RML}.

\subsubsection{One critical point found}

Finding all the critical points is not feasible when the problem  dimension and the complexity  increases. In these cases, it is unlikely that  even the number of critical points will be  known.
Instead of seeking to compute all critical points, we  compute a single critical point using a random starting point for the optimization. 
We assume, without evidence, that the optimization performed this way uniformly samples the critical points, in which case the resulting $\sum_{j=1}^{n_{c}} m^{(j)} w^{(j)}$ provides an unbiased estimator in the same manner as it  is obtained from computing all critical points.

\vspace{\baselineskip}
\begin{algorithm}[H]
\caption{Weighted RML -- sample one critical point}
\SetAlgoLined
\begin{itemize}
	\item[] for $i =1$ to $N_e$
    \begin{itemize}
		\item[] generate samples $m_i'$	 and $\delta_i'$ from $q_{M'\Delta'} (m',\delta') $
			\item[] randomly generate $m_{init}$  (initial guess)
			\item[] solve (\ref{eq:critical}) for $m_{i}$ and set $\delta_i = \delta_i'$
			\item[] compute $w'_{i}$ using \eqref{eq:weight-RML} with $n(m') = 1$
	\end{itemize}
	\item[] assign $w_{i} = w'_{i} /\sum_i w'_{i}$
\end{itemize}\label{Alg3}
\end{algorithm}
\vspace{\baselineskip}

\subsubsection{Local minimizers only}

Due to the high dimension of most realistic applications, solutions of (\ref{eq:critical}) are much more difficult to obtain than local minimizers of $L$, although local maximizers can clearly be easily found as  minimizers of $-L_i$. In general, however, it appears that in high dimensions, minimizers of $L_i$ are far more important than maximizers. When that is the case, a Gauss-Newton approximation of the Jacobian determinant can be made with little loss of accuracy
\begin{equation*}
\begin{split}
J & = |  I +  D (C_M G\trp C_D^{-1}(g(m)-\delta) )  | \\
 & \approx |  I +  C_M G\trp C_D\inv G   |\ .
\end{split} 
\end{equation*}

Thus, in the Darcy flow problem, we will solve for random minimizers instead of solving for random critical points and we will compute a Gauss-Newton approximation of the Jacobian determinant.   The consequence of this approximation is that the distribution of samples will not be exact. In section \ref{sec:quadratic}, we investigate numerically the consequence of sampling only the minimizers.

To obtain the weights, the computation of the Jacobian determinant and the gradient of the misfits with respect to the parameter are necessary. For the first two examples, it is easy to compute the weights. For the large-scale flow problem, we use the Gaussian-Newton method in the hIPPYlib \cite{villa:21} to get a low-rank approximation of  the Jacobian determinant. This can reduce the cost.

\section{Numerical examples}

In this section, we present sampling results using test cases of increasing size and complexity. We first demonstrate the methodology using a toy problem that has been previously used by \cite{wang:18}. It is small enough that the true posterior distribution is easily derived. We focus on the parameter values that are difficult to sample correctly. This example shows the difficulty with using only the minimizers and not accounting for the limited range in the solutions.

A second simple example is the ``banana-shaped'' distribution from \cite{haario:99}. It has been used fairly often to test adaptive forms of MCMC \cite{haario:01,vrugt:06,roberts:09,martino:15,duncan:16}. It is simple enough  that computing the Jacobian determinant is not a challenge so the focus again is on showing that if we find the roots and compute the weights, the sampling is correct.

\subsection{Bimodal posterior pdf}\label{sec:quadratic}

The first example has been previously discussed in  \cite{wang:18} where  it was used to demonstrate properties of the randomized maximum a posteriori sampling algorithm. One of their test problems required sampling from the distribution
\begin{equation*}
\pi_M(m) \propto  \exp \left( - \frac{1}{2} (m-0.8)^2  - \frac{1}{2  \sigma_d^2} (m^2 - 1)^2 \right).
\end{equation*}
Although \cite{wang:18} used three different values of $\sigma_d$, we only show results for the most difficult value,  $\sigma_d = 0.5$. For larger and smaller values, the posteriori distribution is more easily modeled as a mixture of Gaussians and is therefore easier to sample.

\begin{figure}[htbp!]
\centering
\subfloat[Solutions of $\nabla L_i(m) = 0$. Color indicates importance weight.]{\label{fig:quadratic_critical_a}
\includegraphics[width=0.46\textwidth]{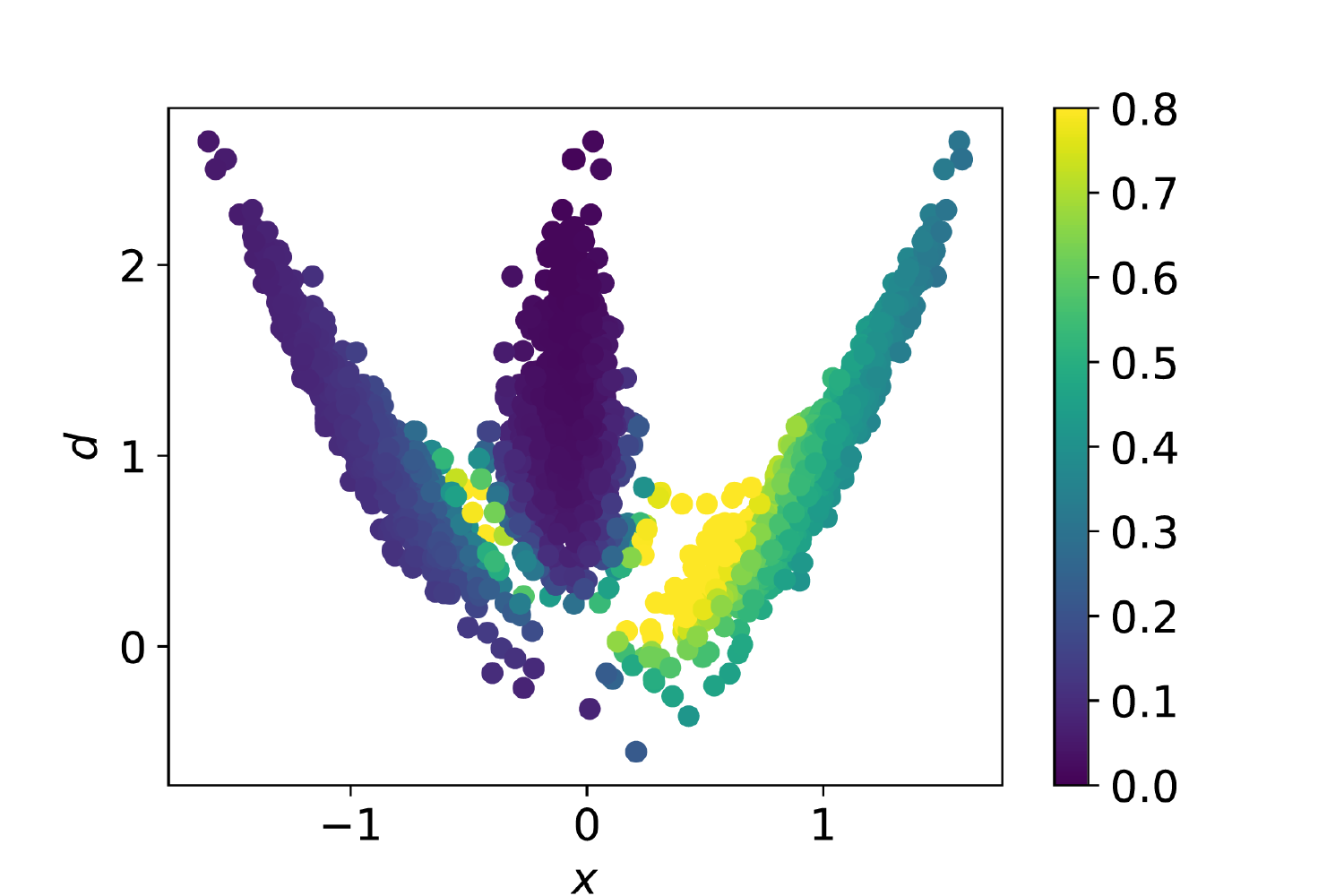}}
~
\subfloat[Distribution of weights from critical points.]{\label{fig:quadratic_critical_b}
 \includegraphics[width=0.46\textwidth]{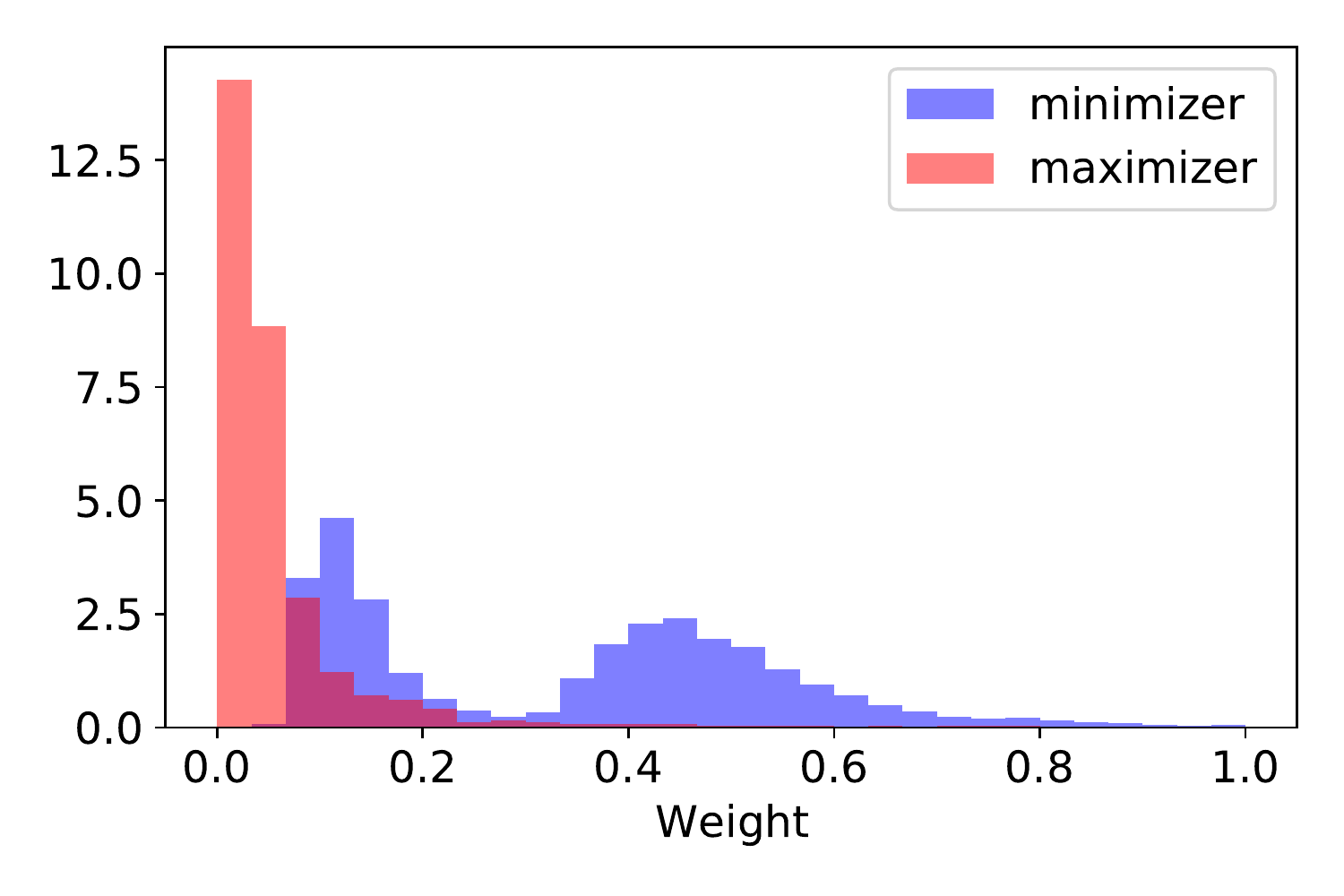}}
\caption{Critical points and importance weights for the quadratic observation operator.}
\label{fig:quadratic_critical}
\end{figure}

In our approach, approximate samples from the posteriori distribution are obtained by solving  for the critical points of
\begin{equation}
L_i(m) = \frac{1}{2} (m-m_i')^2  + \frac{1}{2 \sigma_d^2} (m^2 - \delta_i' )^2
\label{eq:quadratic_objective}
\end{equation}

\[ m_i' \sim N(0.8,1) \qquad \delta_i' \sim N(1, 0.25). \]
Because the objective function in this case is a polynomial, it is straightforward to obtain all real roots of (\ref{eq:critical}). For most choices of $(m',\delta')$ there are three real roots -- $N_s = 10,000$ samples of $(m',\delta')$ from the prior generated $N_e = 25,046$ pairs of $(m,\delta)$. 
The locations of the roots are shown in 
Figure~\ref{fig:quadratic_critical_a}. The set of points in the center of the plot correspond to maximizers of  \eqref{eq:quadratic_objective}. The points on the right side correspond to the global minimum, and the points on the left correspond to the local minimum. The colors show unnormalized importance weights for each sample. The maximizers are generally given small weights (Fig.~\ref{fig:quadratic_critical_b}), although  a small number of maximizers have weights that are similar to the weights of points near the local minimum.

\begin{figure}[htbp!]
\centering
\subfloat[$N_s = 200$.]{\label{fig:quadratic_three_ensemble_sizes_a}
\includegraphics[width=0.32\textwidth]{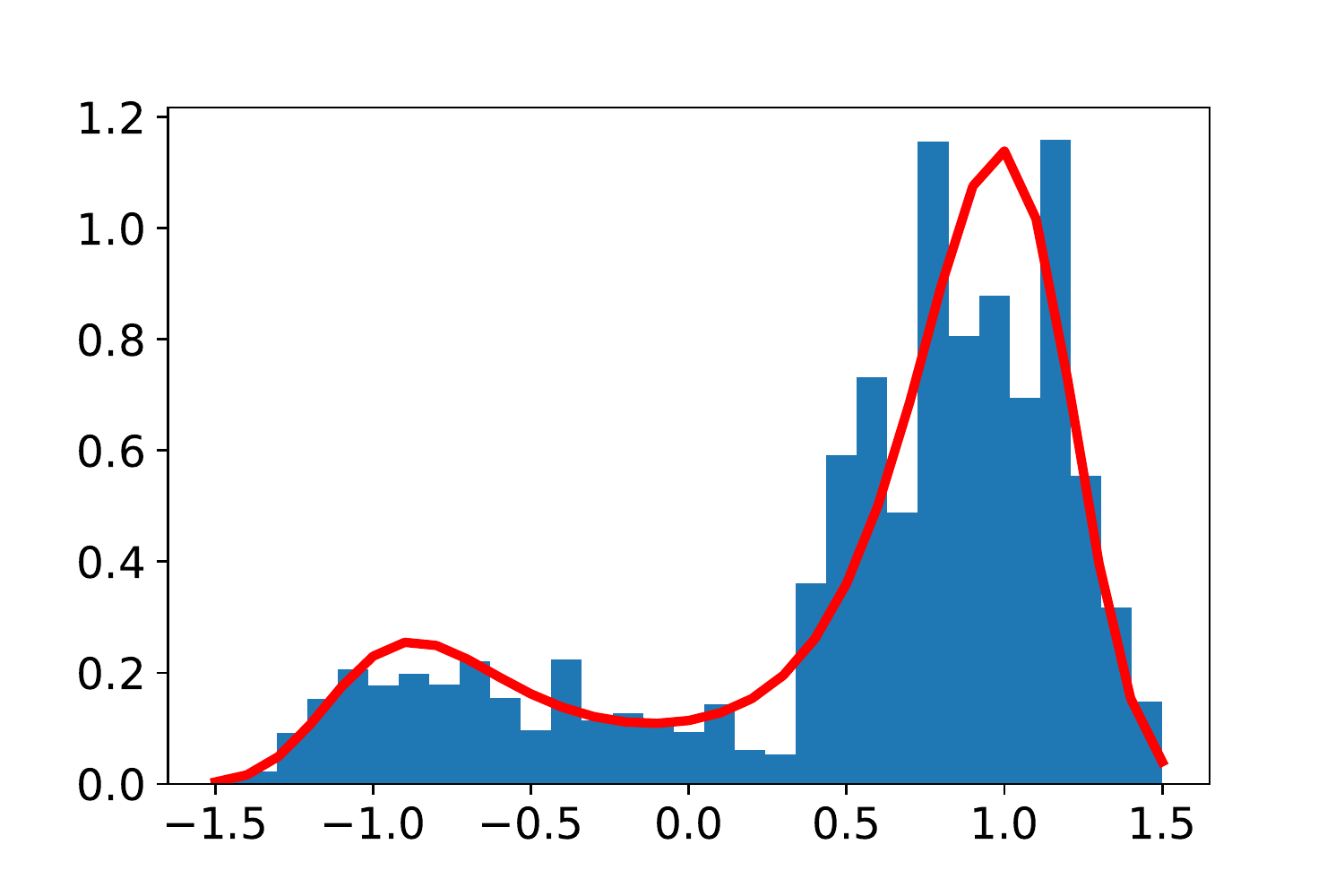}}
    ~
 \subfloat[$N_s = 1000$.]{\label{fig:quadratic_three_ensemble_sizes_b}
 \includegraphics[width=0.32\textwidth]{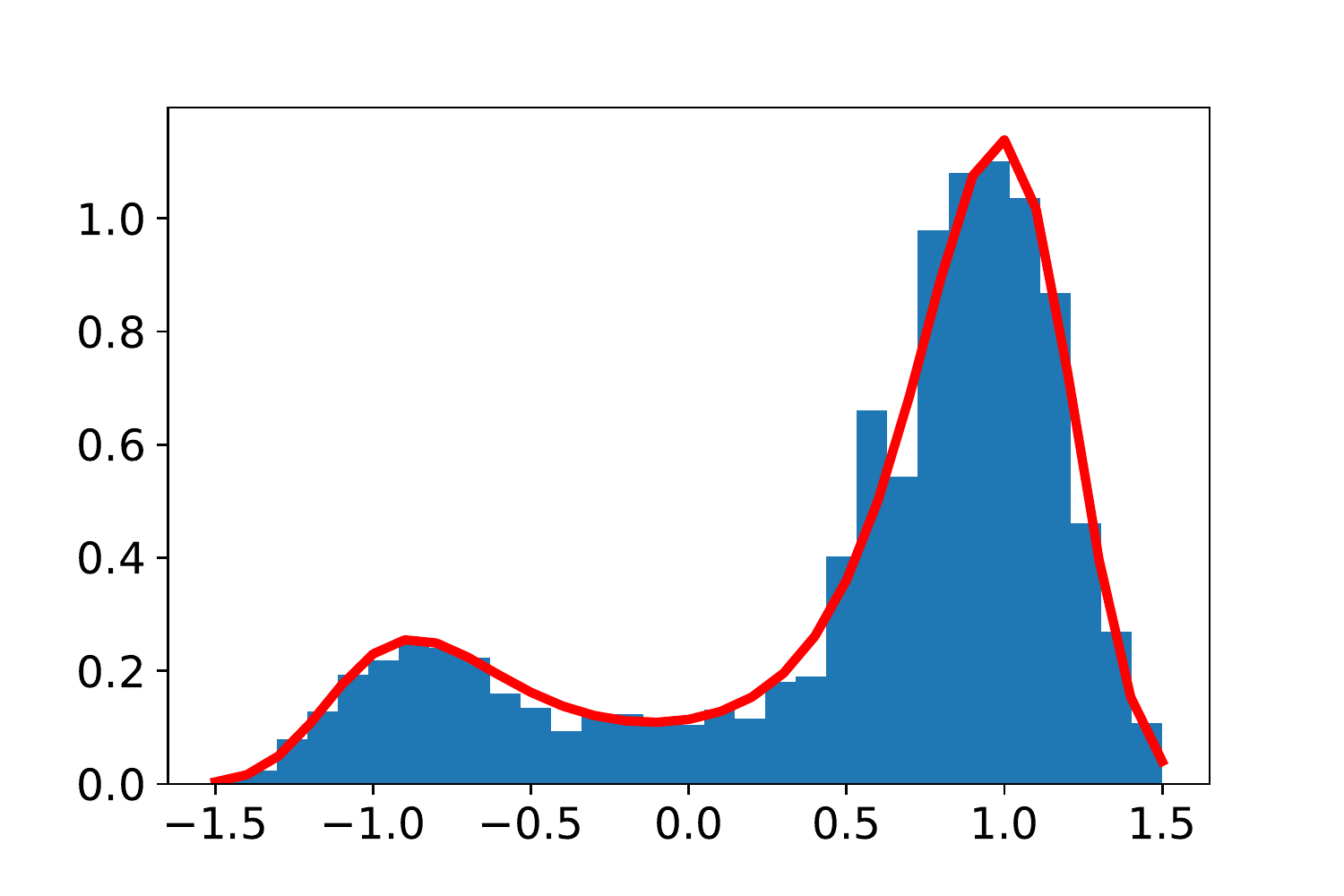}}
    ~
\subfloat[$N_s = 5000$.]{\label{fig:quadratic_three_ensemble_sizes_c}
\includegraphics[width=0.32\textwidth]{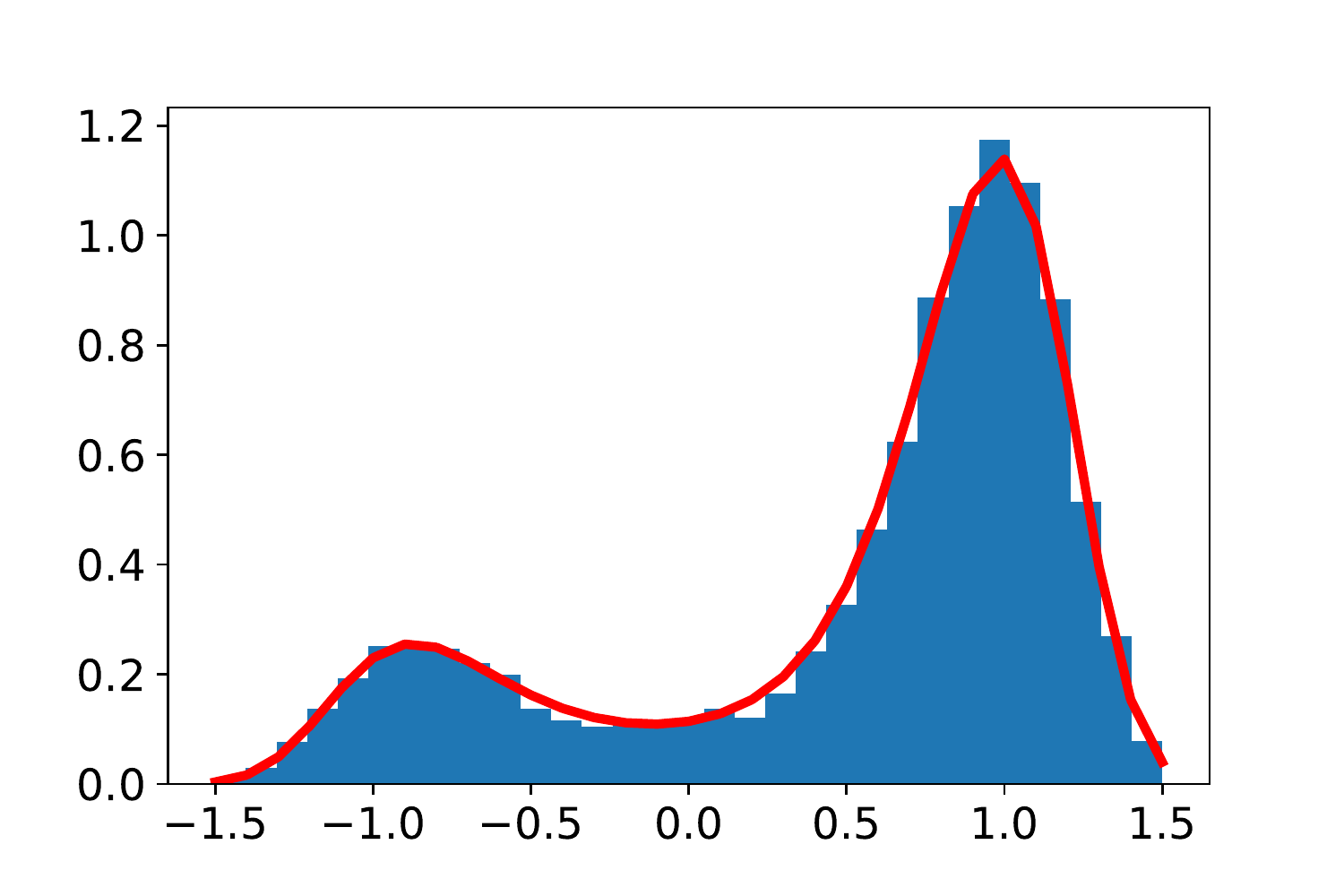}}
\caption{Weighted sampling approximations to the true posterior distribution.}
\label{fig:quadratic_three_ensemble_sizes}
\end{figure}

Because the samples are generated independently (or in groups of three), the quality of the weighted sampling approximation to the target distribution is limited only by sampling error -- larger samples provide better approximations. Figure~\ref{fig:quadratic_three_ensemble_sizes} shows results for three different sample sizes $N_s$ ($200$, $1000$, $5000$).  The number of weighted samples is larger than the ensemble size because a single sample from the prior usually results in three critical points.

\begin{figure}[htbp!]
\begin{tabular}{cccc}
  &  no importance weights  &  GN approx  &  weighted \\
  \raisebox{1ex}{\rotatebox{90}{all critical pts}}
 & \begin{overpic}[width=0.3\textwidth]{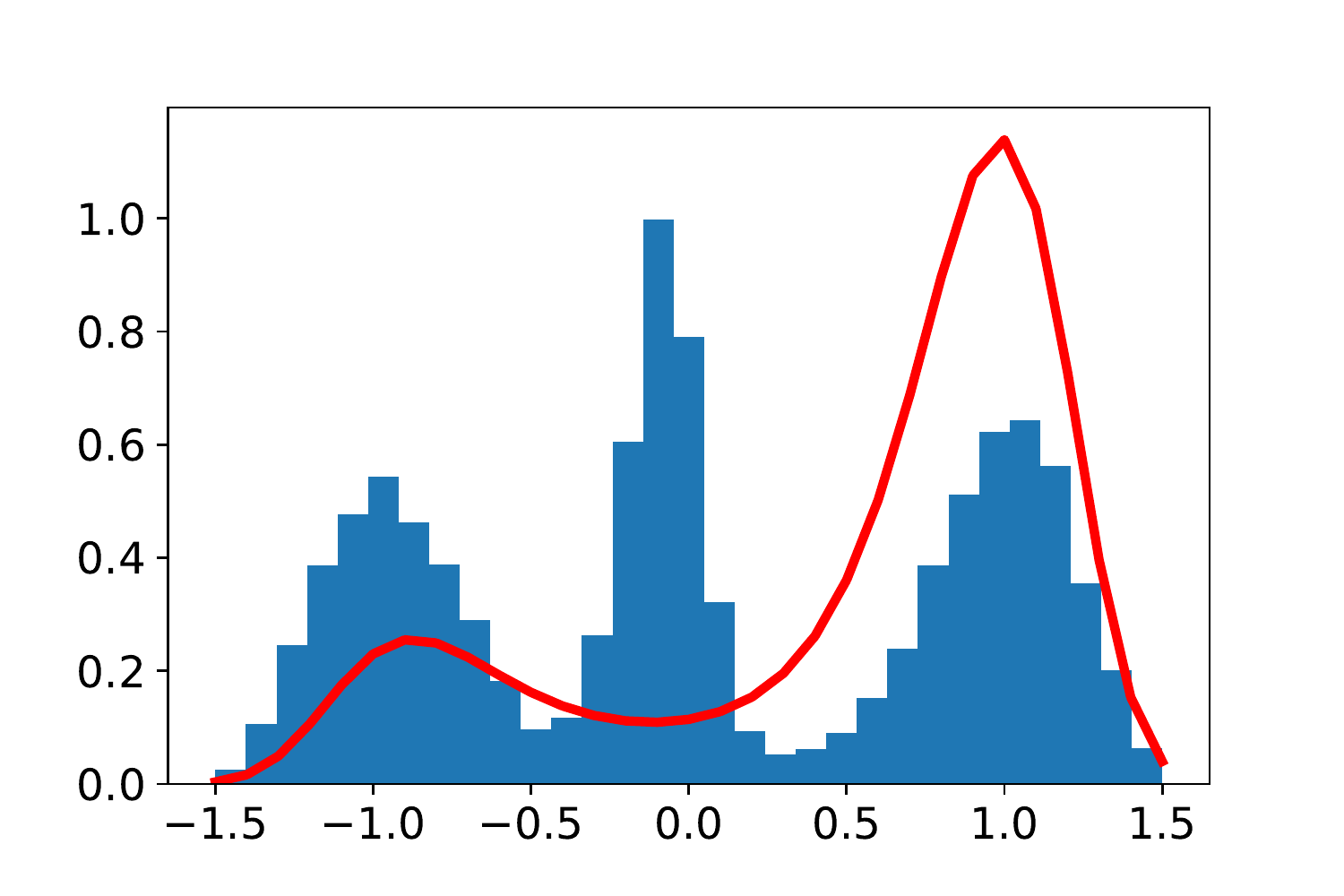}
 \put(15,50){(a)}
 \end{overpic}
& \begin{overpic}[width=0.3\textwidth]{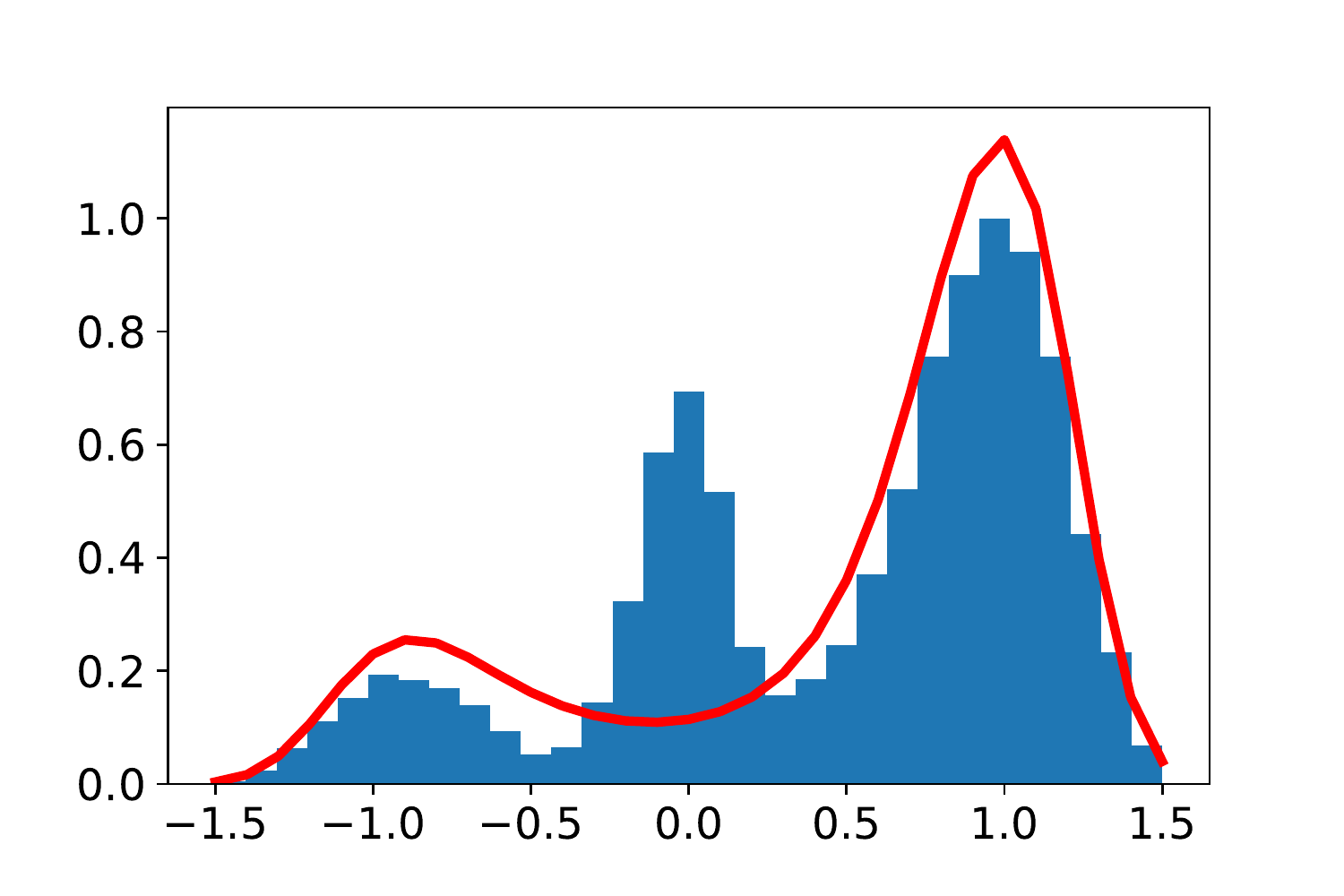}
 \put(15,50){(b)}
 \end{overpic}
& \begin{overpic}[width=0.3\textwidth]{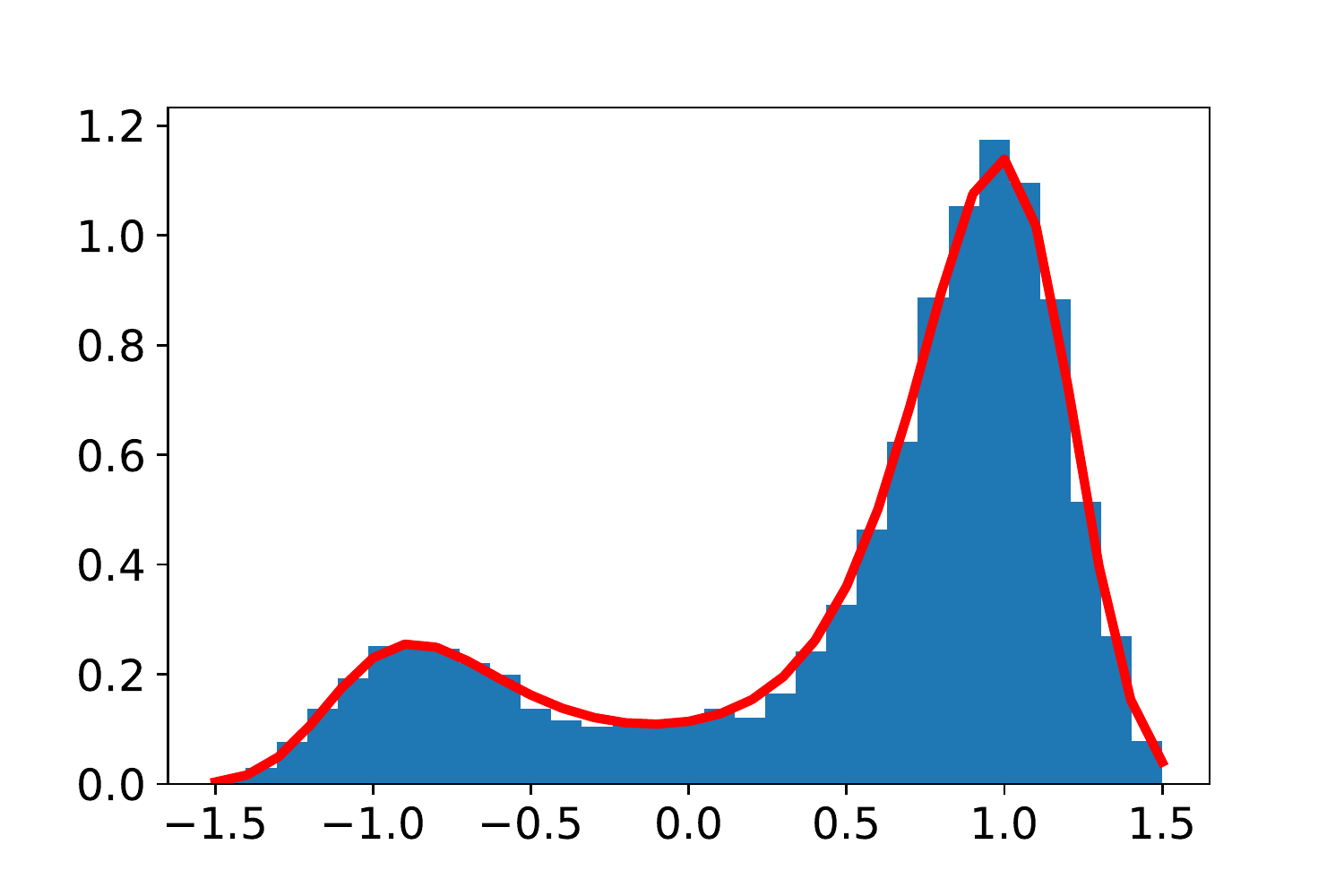}
 \put(15,50){(c)}
 \end{overpic}\\
\raisebox{1ex}{\rotatebox{90}{minimizers}}
 &
\begin{overpic}[width=0.3\textwidth]{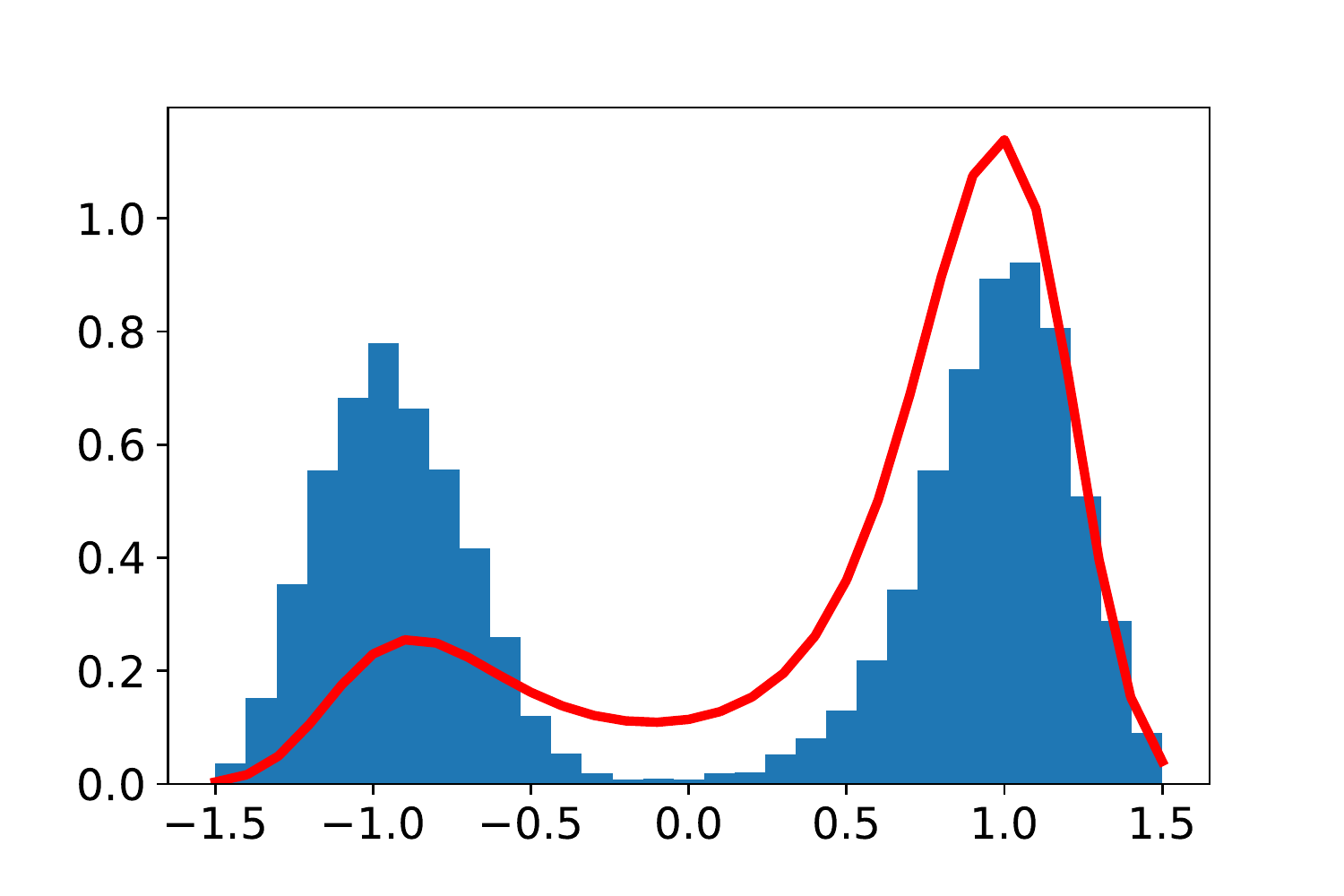}
 \put(15,50){(d)}
 \end{overpic}
& \begin{overpic}[width=0.3\textwidth]{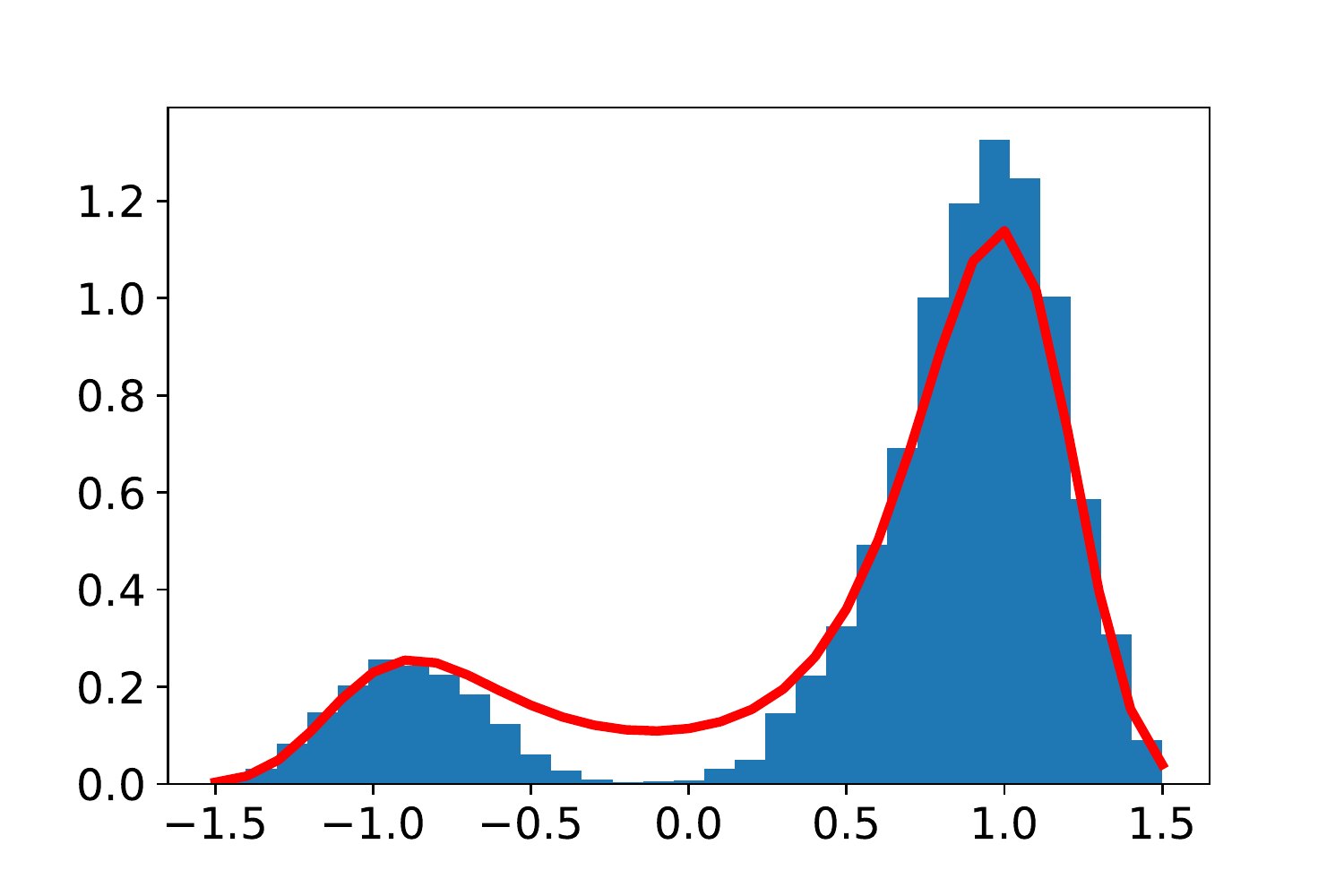}
 \put(15,50){(e)}
 \end{overpic}
& \begin{overpic}[width=0.3\textwidth]{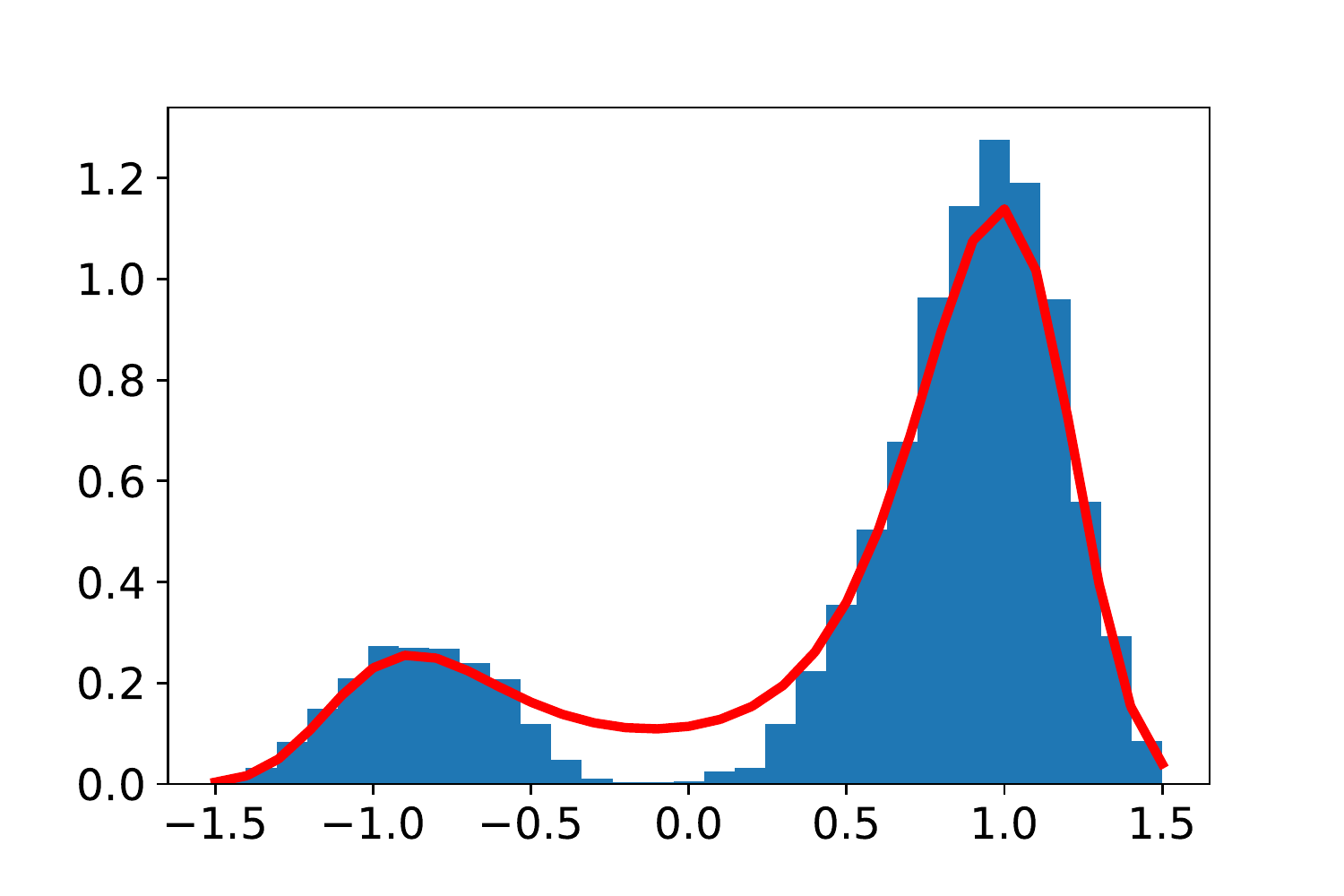}
 \put(15,50){(f)}
 \end{overpic}
\end{tabular}
\caption{Compare distributions from approximations to sampling based on computation of critical points of cost function.}\label{fig:quadratic_three_approx_approaches}
\end{figure}

For this problem, computing all real roots of (\ref{eq:critical}) and computing importance weights for each root is trivial. For large high-dimensional problems, finding multiple roots and computing the weights will be challenging. Here we examine the consequence of three realistic approximations to correct sampling: (1) identifying only the minimizers of the cost function, (2) using a Gauss-Newton (GN) approximation of the Jacobian of the transformation and (3) neglecting importance weights altogether. Fig.~\ref{fig:quadratic_three_approx_approaches}c shows that correct sampling of the target distribution is obtained when all critical points are included and the weights are computed accurately. If the importance samples are neglected (Fig.~\ref{fig:quadratic_three_approx_approaches}a) or if the GN approximation of the Jacobian is used (Fig.~\ref{fig:quadratic_three_approx_approaches}b), the distribution of samples is badly distorted. When it is not possible to compute the Jacobian accurately in high dimensions, it appears to be advisable to only compute the minimizers. The distribution of samples obtained using the GN approximation applied to minimizers (Fig.~\ref{fig:quadratic_three_approx_approaches}e) is nearly as good as the results with correct weights, and far better than results with  no importance weighting.

\subsection{Banana-shaped posterior pdf}

The second numerical example is the widely used
``banana-shaped'' target density
initially presented in \cite{haario:99}, but extended to higher dimensions by \cite{roberts:09},
\begin{equation}
\pi_M(m) \propto  \exp \left[ - \frac{1}{2 \sigma_m^2 } \left(  m_1^2 + m_2^2 + \cdots +  m_{N_m}^2 \right) \right]  \exp \left[ - \frac{1}{2 \sigma_d^2 }  (4 - 10 m_1 -m_2^2 )^2   \right] 
\label{eq:roberts09}
\end{equation}
with   $\sigma_d = 4$ and $\sigma_m = 5$.   In our numerical experiment, we use parameters values from \cite{martino:15}, but increased the dimension of $m$ to 4.
The  first term of \eqref{eq:roberts09} is identified as the  Gaussian prior  with model covariance $C_M = I$ and the second term as the log-likelihood, with
$$
g(m) = 10 m_1 + m_2^2.
$$

\begin{figure}[htbp!]
\centering
\subfloat[$m_1$-$m_2$ plane.]{\label{fig:banana_minimizers_a}
\includegraphics[width=0.30\textwidth]{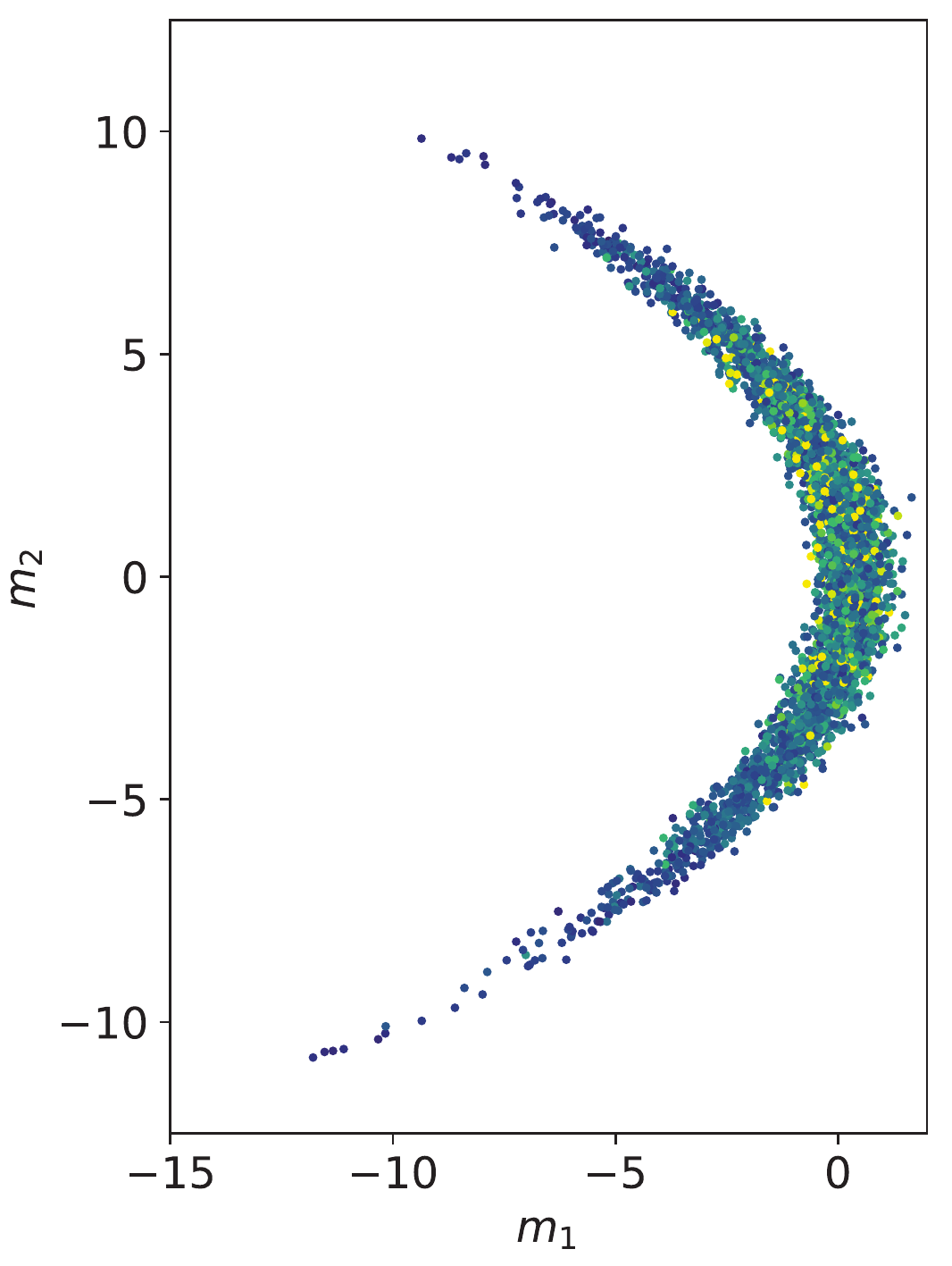}}
    ~
\subfloat[$m_1$-$m_3$ plane.]{\label{fig:banana_minimizers_b}
\includegraphics[width=0.26\textwidth]{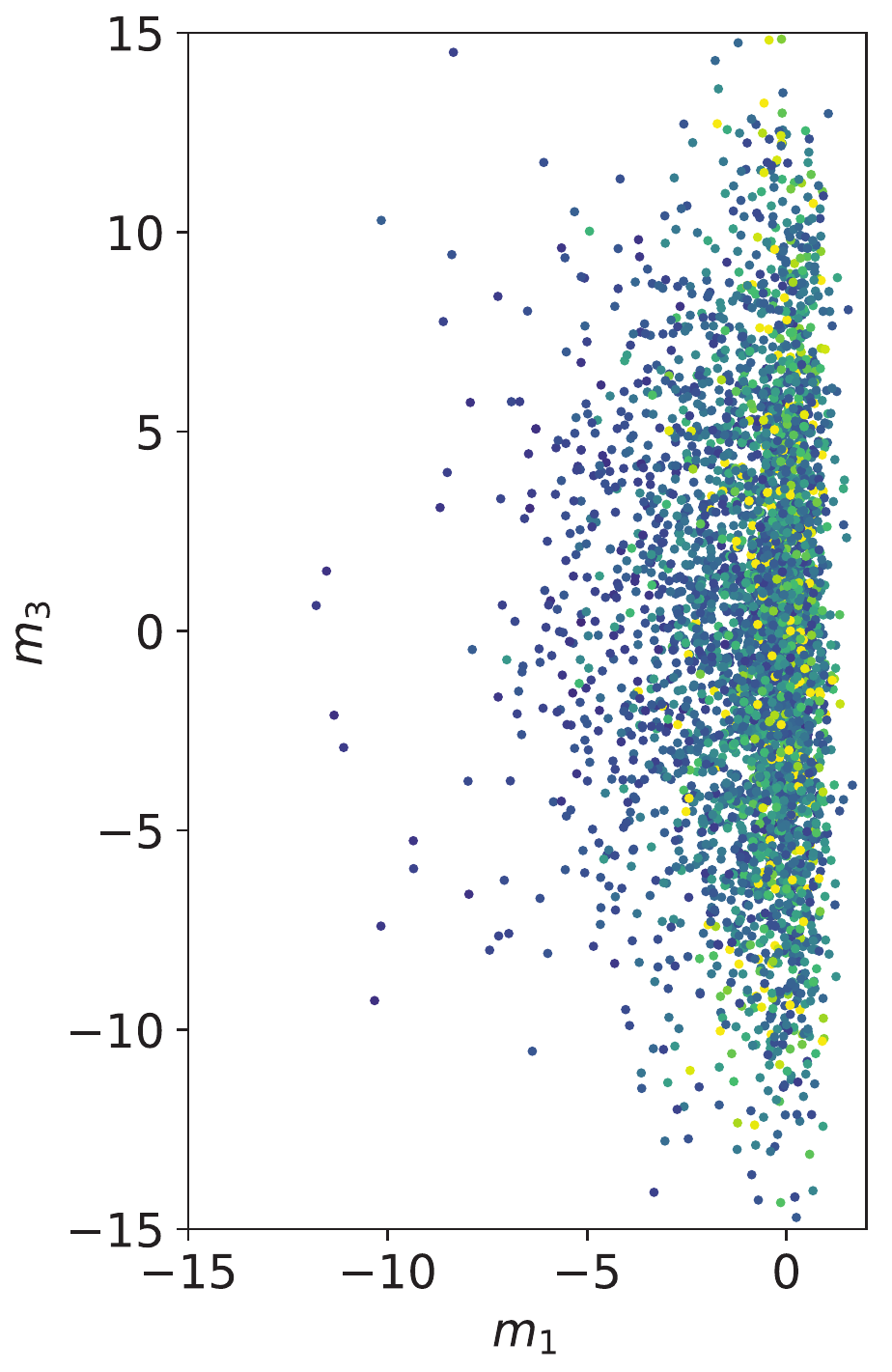} }
    ~
\subfloat[$m_1$-$\delta$ plane.]{\label{fig:banana_minimizers_c}
\includegraphics[width=0.36\textwidth]{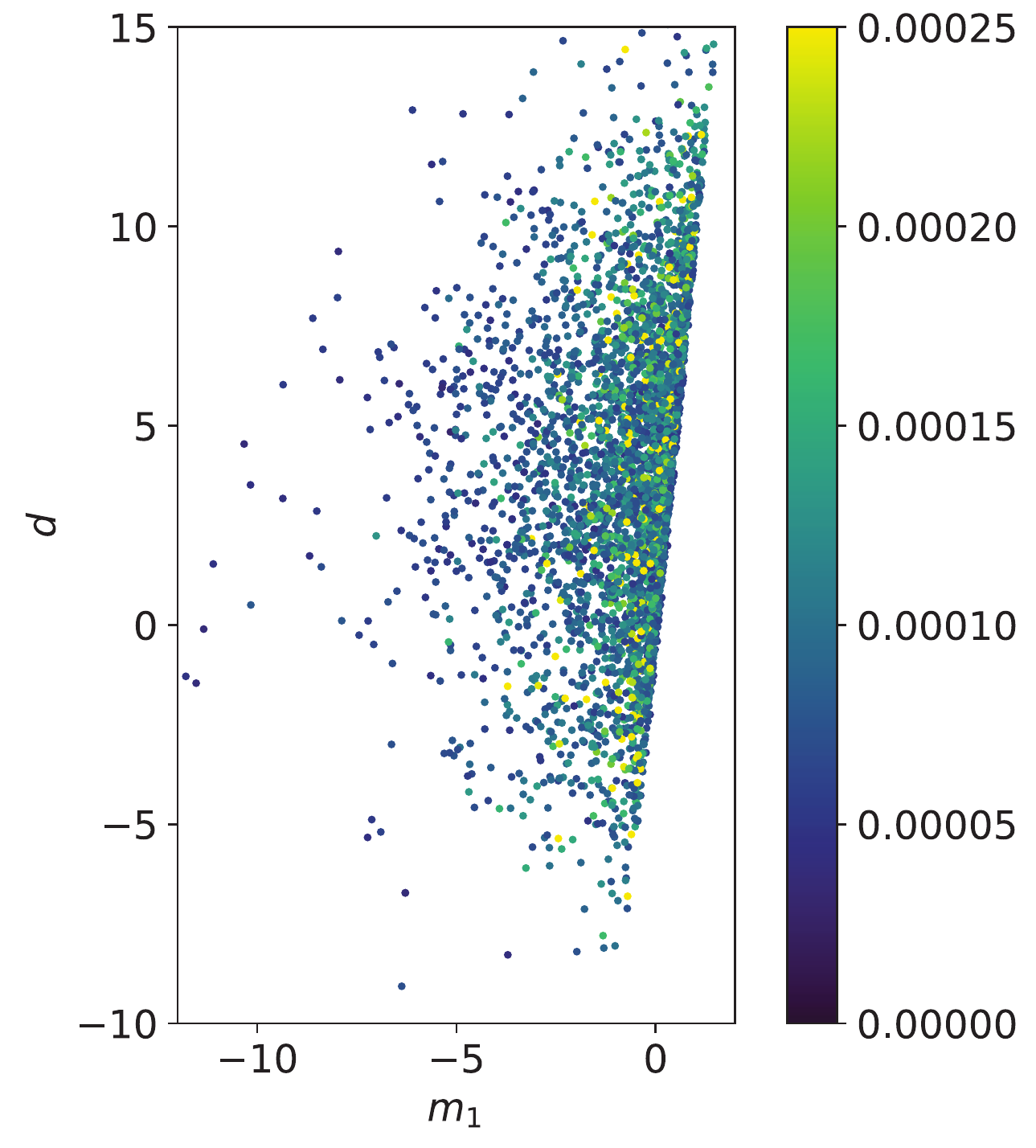} }
\caption{Critical points of the banana-shaped objective function. Color indicates importance weight on the samples. The same color scale is applied to all subplots.}
\label{fig:banana_minimizers}
\end{figure}

\begin{figure}[htbp!]
\centering
\subfloat[Distribution of importance weights for the minimizers of the log posterior.]{\label{fig:banana_contour_a}
\includegraphics[width=0.56\textwidth]{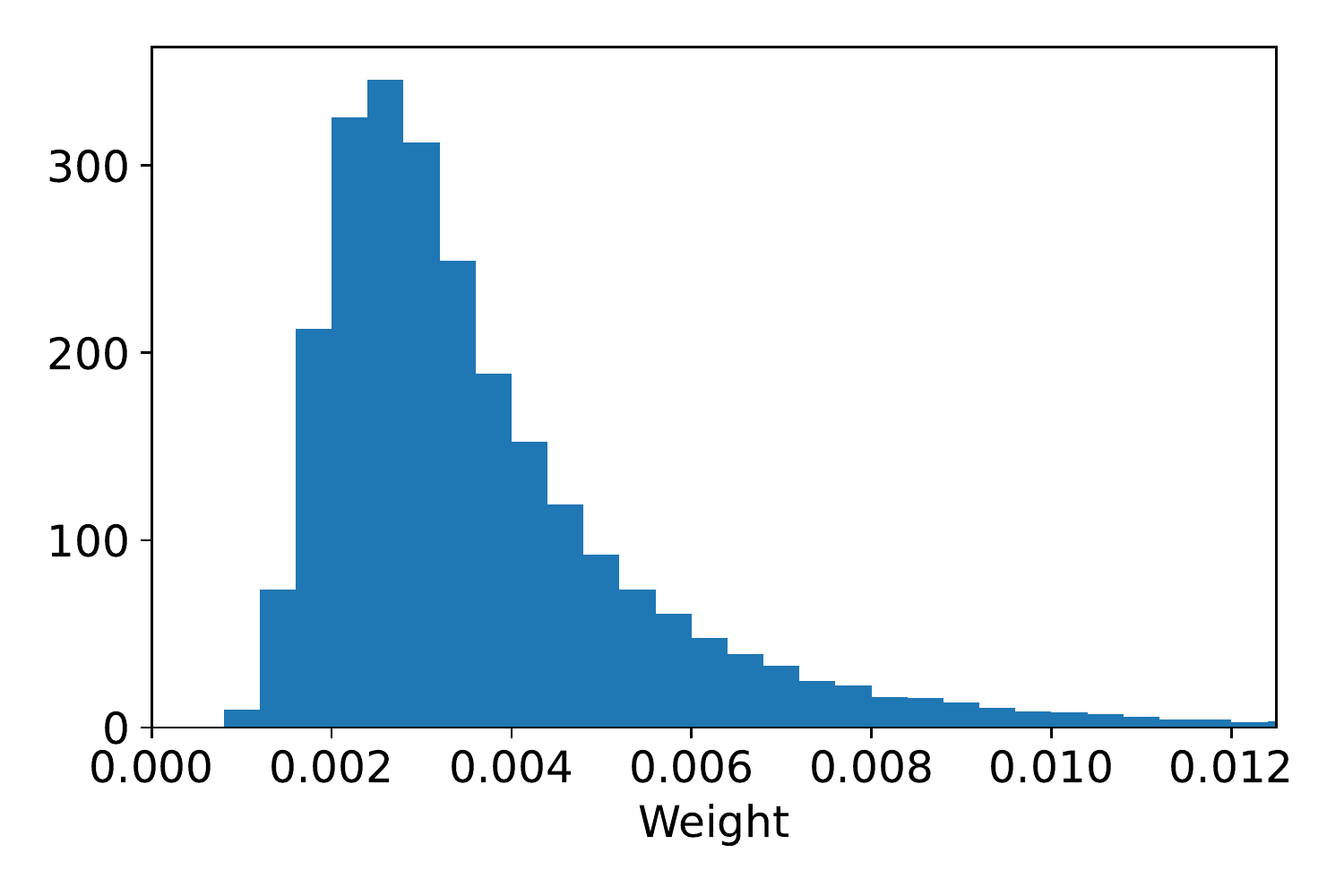}}
~
\subfloat[Compare true distribution with estimation from weighted minimizers ($m_1$-$m_2$ plane).]{\label{fig:banana_contour_b}
\includegraphics[width=0.36\textwidth]{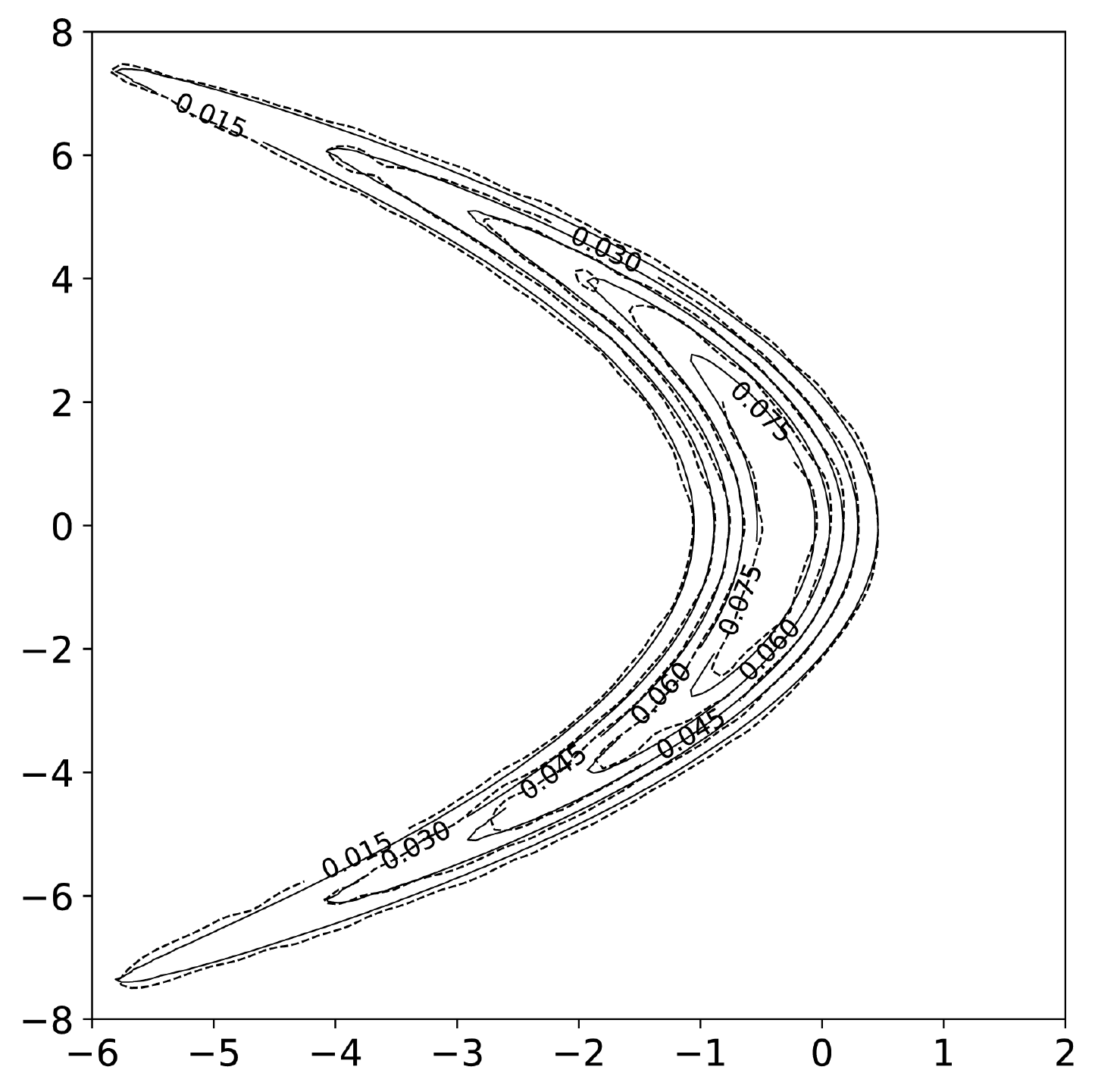} }
\caption{Minimizers of the banana-shaped objective function.}
\label{fig:banana_contour}
\end{figure}

Because of the curved shape of the objective function (Fig.~\ref{fig:banana_contour_b}), accurate computation of minimizers of $L_i$ was relatively difficult.  Three projections of  the first 3000 approximate samples obtained using the Broyden-Fletcher-Goldfarb-Shanno algorithm for minimization  are shown in Fig.~\ref{fig:banana_minimizers}. Each minimization  was initiated at the point $m_i'$.

For this problem, which has a single critical point,  the empirical distribution obtained from minimization appears to be relatively good.
The true conditional distribution distribution  $\pi(m_1,m_2|m_3=0,m_4=0)$ is compared in Fig.~\ref{fig:banana_contour_b}, with a kernel estimate of the empirical marginal distribution for $m_1, m_2$ (dashed contours) obtained using 50,000 minimizations. Also, based on the distribution of unnormalized importance weights  (Fig.~\ref{fig:banana_contour_a}), it appears that the weights are not dominated by a few large values, which is confirmed by a high effective sampling efficiency, $N\sbr{eff}/N_e = 44796/50000 \approx 0.9$, based on Kong's estimator \eqref{eq:NEff},
\begin{equation}\label{eq:NEff}
N_{\text{Eff}}=\frac{1}{\sum_{k=1}^{N_e}w_k^2},
\end{equation}
where $\sum_{k=1}^{N_e}w_k=1$

\subsection{Darcy flow example} \label{sec:flow}

\begin{figure}[htbp]
\begin{tabular}{cc}
\includegraphics[width=0.4\textwidth]{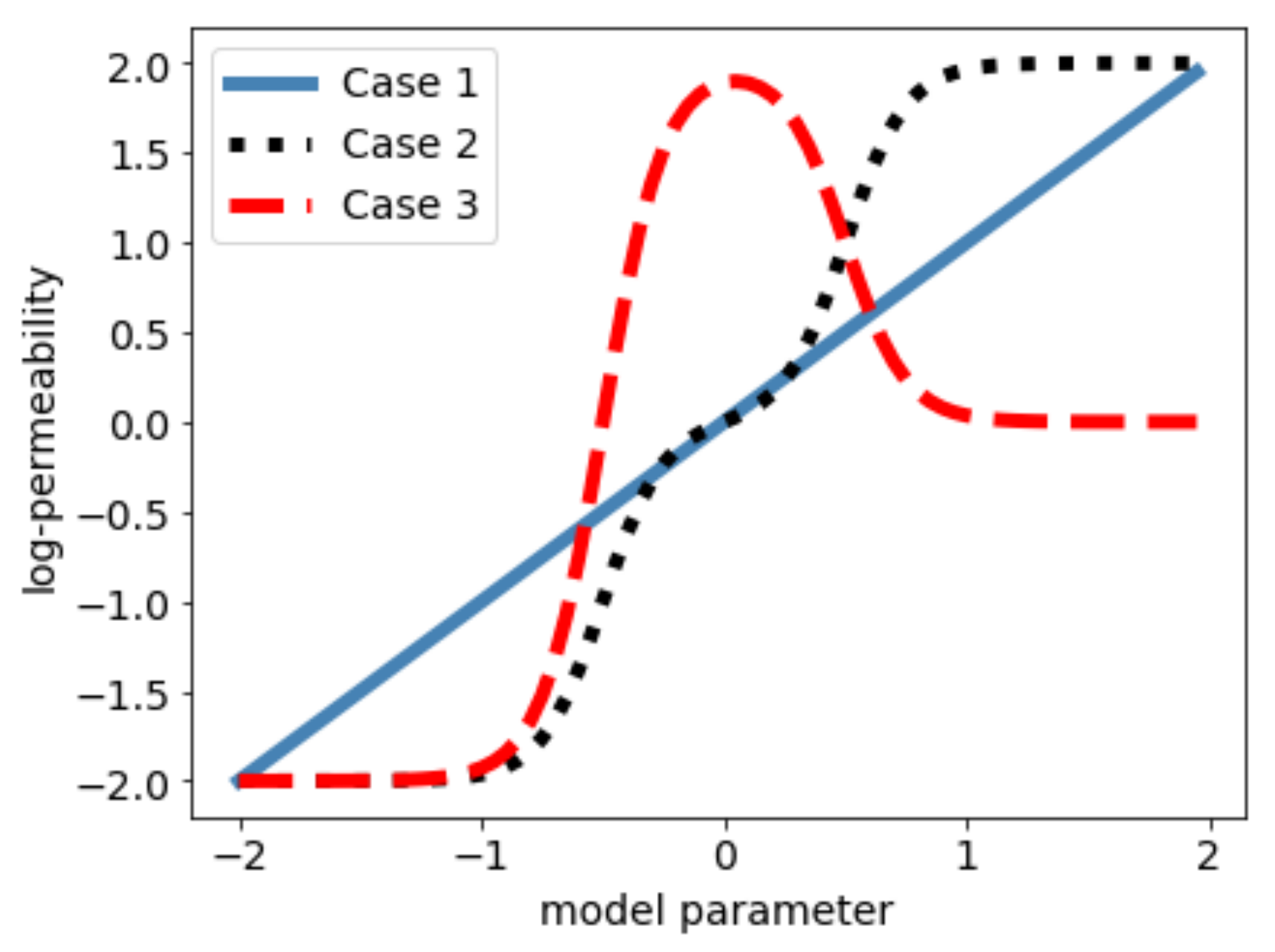}
&
\includegraphics[width=0.45\textwidth]{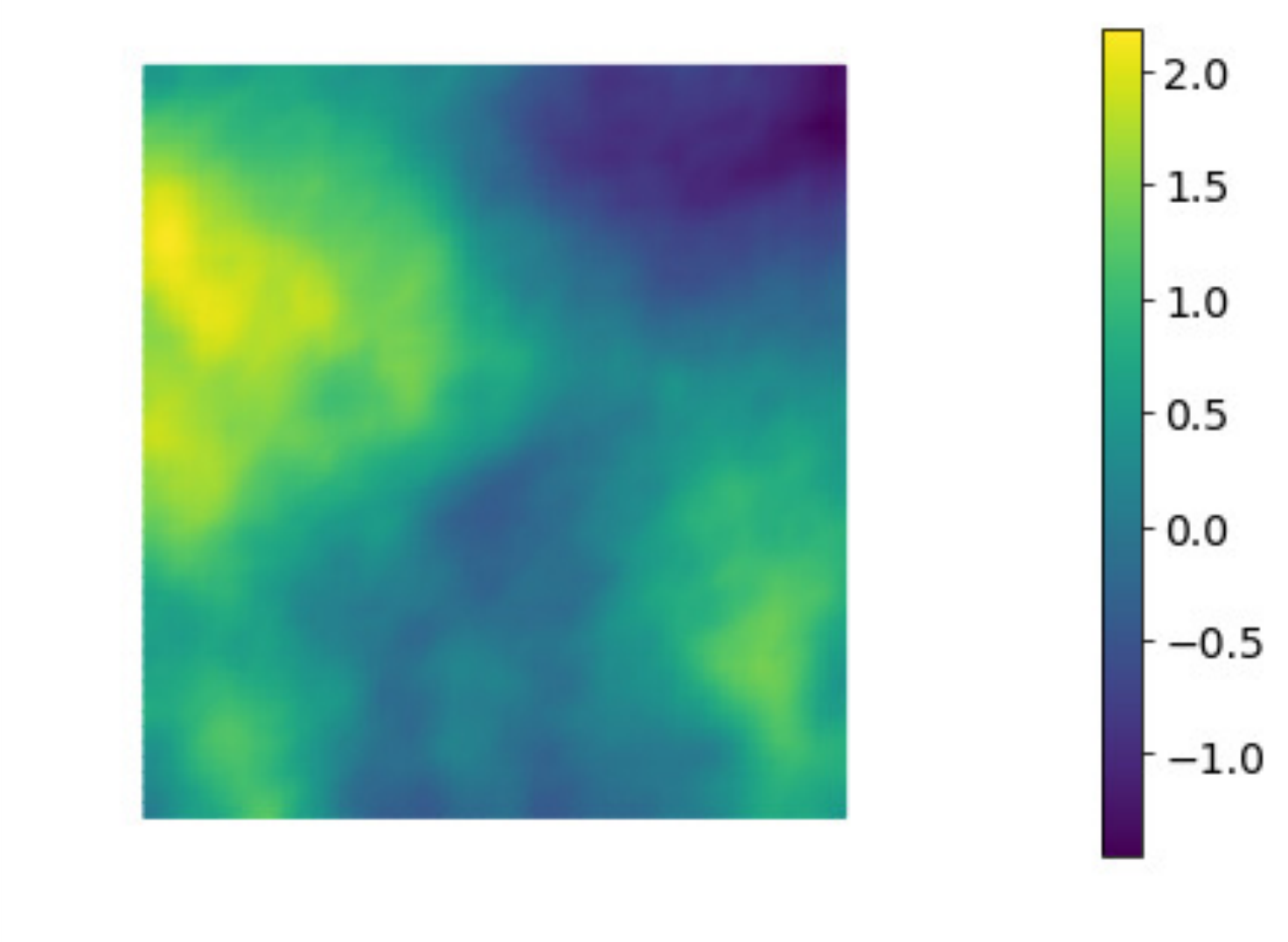} \\
(a) Three transformations & (b) Case 1 \\[1em]
\includegraphics[width=0.45\textwidth]{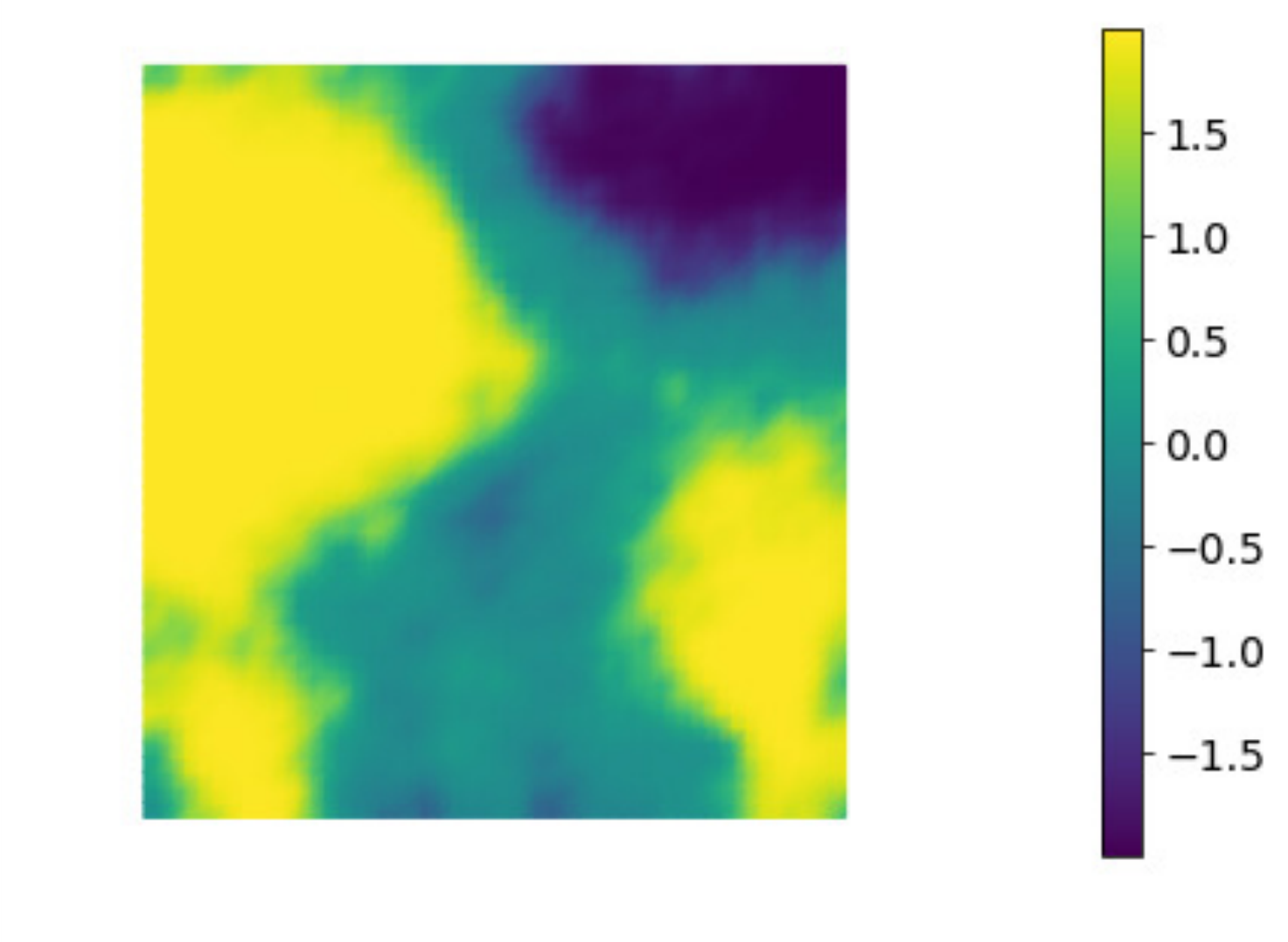}
&
\includegraphics[width=0.45\textwidth]{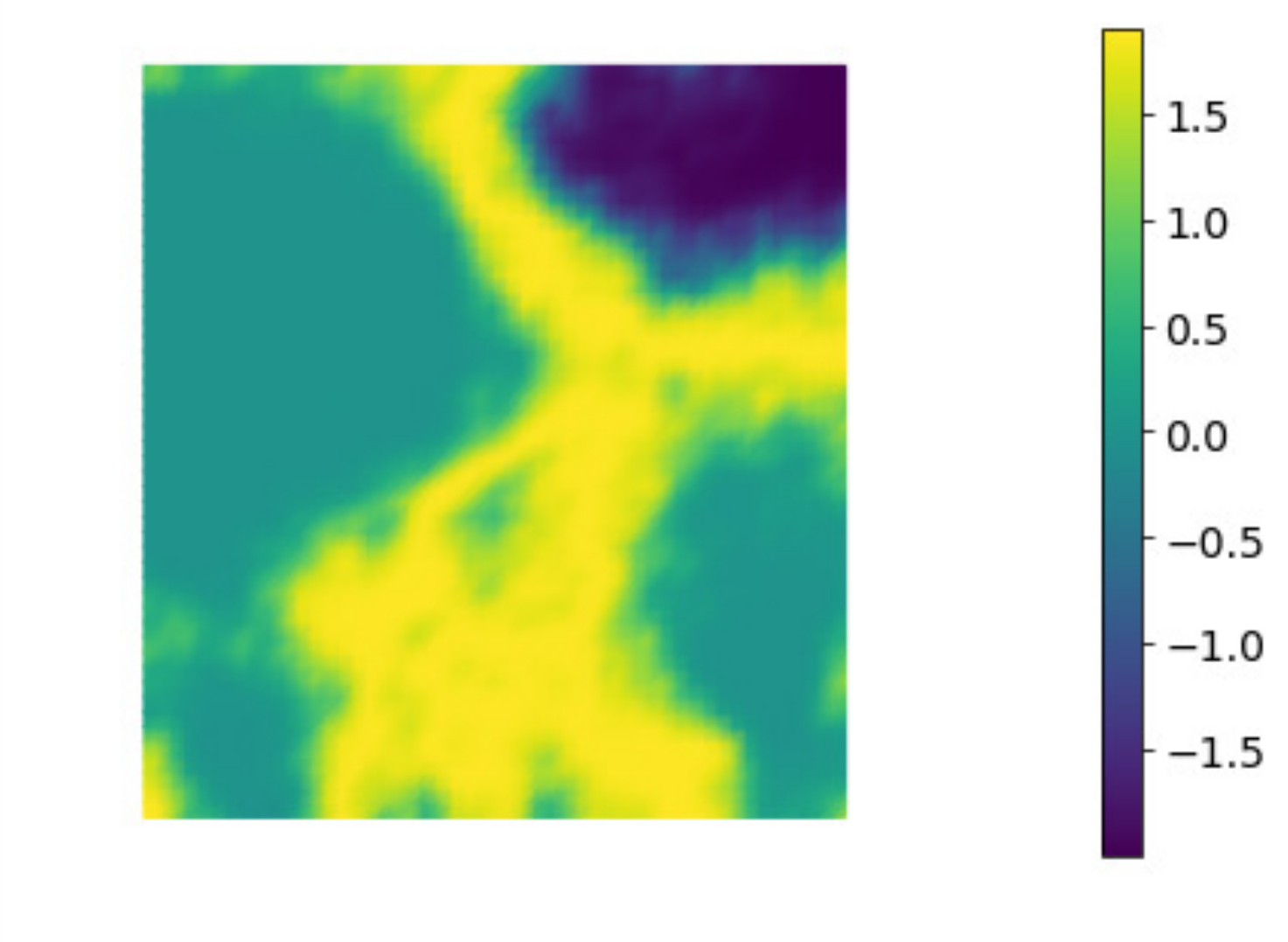} \\
(c) Case 2   & (d) Case 3
\end{tabular}
  \caption{The true log-permeability fields for Cases 1, 2  and 3. All cases use the same true latent variable field.} \label{fig:True}
\end{figure}

For the first two examples, there would be no advantage in using RML for sampling -- MCMC  with a carefully chosen transition kernel would probably be a better alternative in either case. The advantage for RML occurs in large high-dimensional problems for which traditional methods are impractical.
In this section, we investigate the ability to quantify uncertainty in the permeability field $\kappa (x)$ from spatially distributed observation of steady-state pressure.
The pressure $u(x)$ in this example is governed by the equation
\[
-\nabla\cdot(\kappa (x)\nabla u(x)) = 0 \quad \text{in}\quad \Omega=[0,1]\times[0,1]
\]
with the mixed boundary conditions
\[
\left\{
\begin{aligned}
\nabla u\cdot {\mathbf n}&=0 \quad\text{on} \quad\Gamma_{N_{1}}=0\times[0,1]\cup1\times[0,1]\\
\nabla u\cdot {\mathbf n}&=v(x) \quad\text{on} \quad \Gamma_{N_{2}}=[0,1]\times 1\\
u(x) &= 0 \quad\text{on}\quad \Gamma_{D}=[0,1]\times 0.\\
\end{aligned}
\right.
\]
As permeability is a positive quantity, it cannot be modelled as a Gaussian random variable. Here, we evaluate sampling with three possible prior distributions for permeability. In all cases, we define a latent variable $m(x)$ that is multivariate Gaussian, with a prior given by \eqref{eq:prior}.
We take $\alpha=0.12$ and $\gamma=1.12$ which results in a correlation length of approximately 2 and a variance of 1. 
In the first case (Case 1), permeability is modeled as being log-normally distributed, i.e., $\kappa (x) = \exp \big(m(x)\big)$, which is a  typical assumption for the  distribution   of permeability within a single rock type  \cite{freeze:75}.
In more complex formations, it is often useful to model permeability as being largely determined by rock type. In that case, permeability might be largely uniform within a rock type, but variable between rock types. We created two soft thresholding transformations to model the distribution of permeability in a formation with three rock types.
In the the first of the distributions (Case 2), the permeability is related to the latent variable through a highly nonlinear, but monotonic transformation
\begin{equation}\label{eq:threshold2}
\kappa (x) = \exp \Big(  \tanh \big(4 m(x) + 2 \big)  +   \tanh \big(4 m(x)-2\big) \Big) .
\end{equation}
In the second distribution (Case 2), the permeability is related to the latent variable through a non-monotonic transformation
\begin{equation}\label{eq:threshold3}
\kappa (x) = \exp \Big( 2 \tanh \big(4 m(x) +2 \big)  +   \tanh \big(2-4m(x)\big)-1\Big) ,
\end{equation}
which gives a permeability field with a low permeability `background' and connected high perm 'channels' as might occur in subsurface rock formations \cite{armstrong:11}.
Figure~\ref{fig:True} shows the transformations and the three synthetic true log-permeability fields that are used to generate observations for Case 1 (top right), Case 2 (lower left) and Case 3 (lower right).

\begin{figure}
  \centering
    \includegraphics[width=0.95\textwidth]{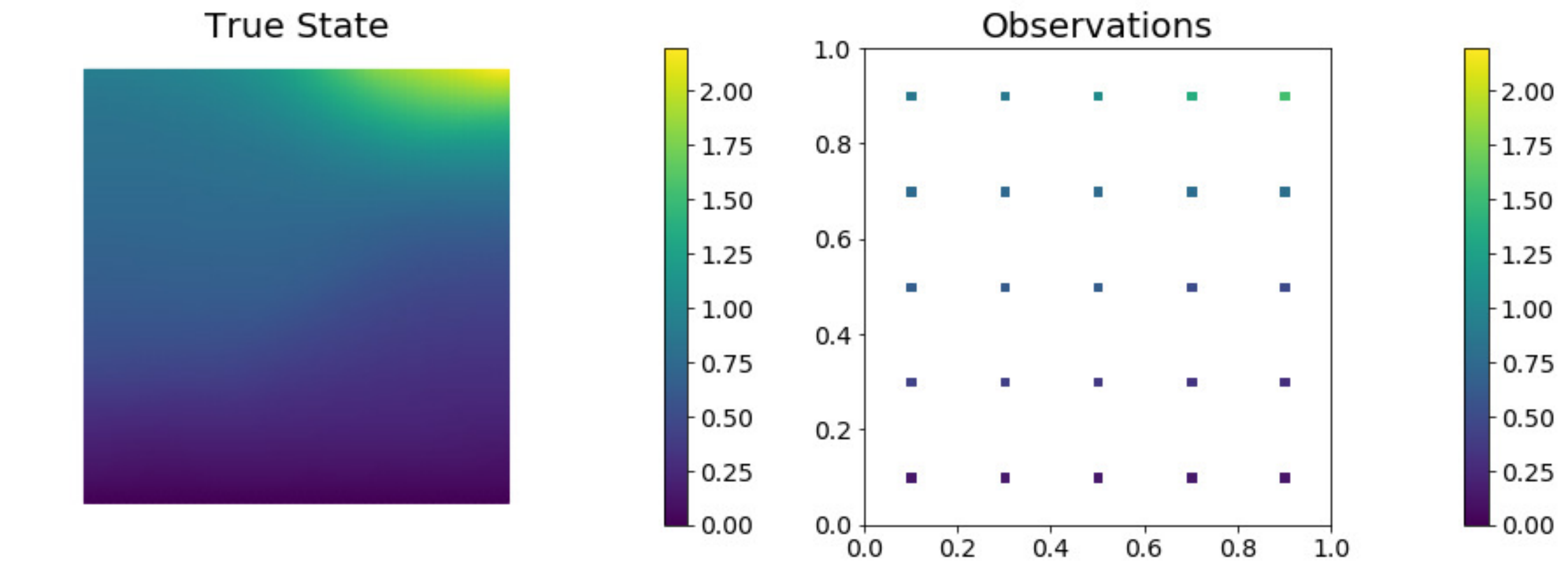}
  \caption{The true state (left) and observation locations (right).}\label{fig:Obs1}
\end{figure}

Figure~\ref{fig:Obs1} (left) shows the true pressure field for Case 1. The pressure distributions for Cases 2 and 3 look similar. For each case, we take 25 pressures as observations. The observation locations are distributed on the uniform $5\times 5$ grid of the domain $[0.1, 0.9]\times [0.1, 0.9]$ as shown in Figure~\ref{fig:Obs1} (right). The noise in the observations is assumed to be Gaussian and independent with standard deviation 0.01.
For the mixed boundary condition, we take $v(x)=2$ for Case 1 and $0.7$ for Case 2. Here the piecewise quadratic finite element is used for the state and adjoint spaces, while piecewise linear finite element is used for the parameter space. The forward model is solved by the finite element method with a uniform $50\times50$ grid for the three cases. Thus the dimension of the discrete state and adjoint space is $10201$ and the dimension of the parameter space is $2601$.

We used the hIPPYlib environment \cite{villa:18,villa:21} for computation of low-rank approximations of eigenvalues of the Hessian.
hIPPYlib
builds on FEniCS \cite{logg:12,langtangen:16} for the discretization of the PDE and uses PETSc \cite{abhyankar:18,zhang:19} for scalable and efficient linear algebra operations and solvers. Minimizers of the objective functions are computed using an inexact Newton-CG solver. We used default parameters for minimization, except that we increased the maximum number of iterations to $300$. The actual average number of iterations required for convergence varied considerably for the three cases.
In Case 1, an average of $24$ iterations were required, Case 2 required $34$ iterations and Case 3 required an average of $73$ iterations.

For the Darcy flow examples, the low-rank approximation of the Jacobian determinant $J$ was used to reduce the cost of computation of weights. The low-rank approximation adopts the Gauss-Newton method in the hIPPYlib \cite{villa:21}. To illustrate the performance of the proposed method, we compare the results obtained by using the stochastic Newton method \cite{martin:12} implemented in hIPPYlib  with that of unweighted RML and weighted RML. Here the stochastic Newton samples are generated from the Gaussian approximation of the posterior. The observation data, prior covariance operator and sample size $N_e$ are the same for the RML and stochastic Newton methods. For convenience, we write the stochastic Newton as SN, unweighted RML as RML and weighted RML as WeRML in the figures of the three cases.

\subsubsection{Case 1: permeability field is log-normal}

As the permeability transformation, $\kappa = \exp{m}$, is monotonic in this example, we might expect the stochastic cost function $L_i$ to have a  single critical point for each sample from the prior. We confirmed this empirically through an investigation in which we generated a single sample from the prior, but $50$ randomly sampled starting points for the minimization.
In the experiments, all $50$ initial starting models converged to the same   model parameters.
In the results that we present, RML sampling was performed with $1000$ samples from the prior and a single random starting point for each minimization. We dropped $15$ out of $1000$ samples for which the minimization routine failed to converge to a sufficiently small value of the gradient norm in $300$ iterations.

The effective sample size computed from \eqref{eq:NEff}, was relatively high for this case:  $N\sbr{eff} \approx 823$ effective samples. The effective sample efficiency, $N\sbr{eff}/N_e = 823/985 = 0.836$. Unlike the toy examples, computing the minimizers and the Jacobian determinant is necessarily approximate in the Darcy flow problem. The expected value for the squared data misfit (with respect to the actual observed values of pressure) is $n_d \sigma_d^2 = 0.0025$, which is somewhat smaller than the mean of the actual squared misfits, $0.0032$.

\begin{figure}[htbp!]
\centering
\subfloat[Case 1 (lognormal permeability)]{\label{fig:Wemi_a}
\includegraphics[width=0.45\textwidth]{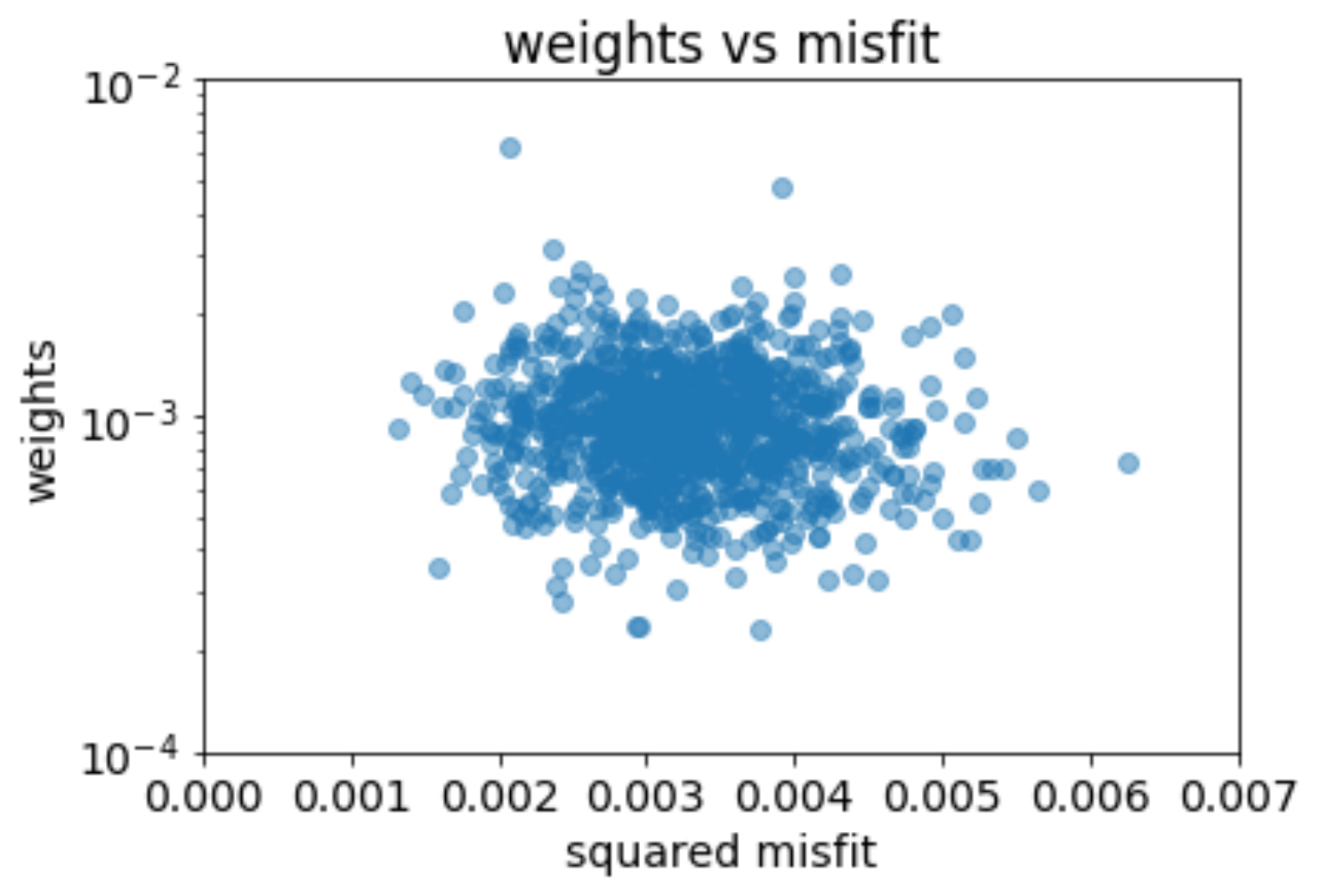}}
~
\subfloat[Case 2 (monotonic log-permeability)]{\label{fig:Wemi_b}
 \includegraphics[width=0.45\textwidth]{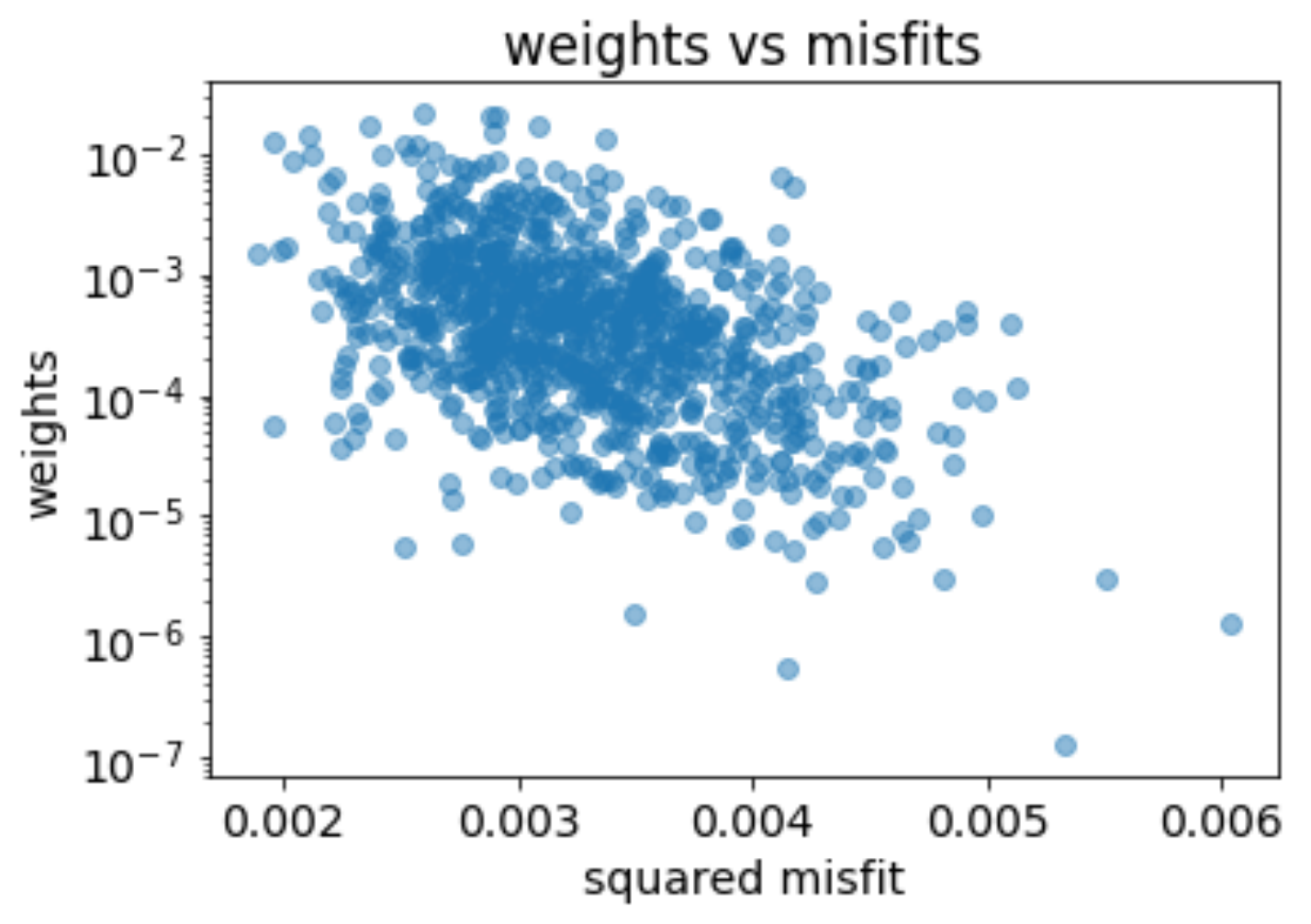}}
\caption{The weights vs misfits for the minimizers. Blue points show computed weights.}
\label{fig:Wemi}
\end{figure}

Figure~\ref{fig:Wemi_a} shows crossplots of the weights vs squared data misfit computed using $V$ as in \eqref{eq:V}.
The  approximation of log-weights is only slightly correlated with the squared data misfit  ($r=-0.09$) so it appears that in this case, as in the quadratic example,  the data mismatch would not serve as a viable surrogate for weighting of samples.
In addition to data mismatch, the weights are affected by the  nonlinearity of the problem -- either through $|V|$ or through the term  $(g(m) - Gm)$, which occurs in $\eta$.

\begin{figure}[htbp]
  \centering
  \includegraphics[width=0.95\textwidth]{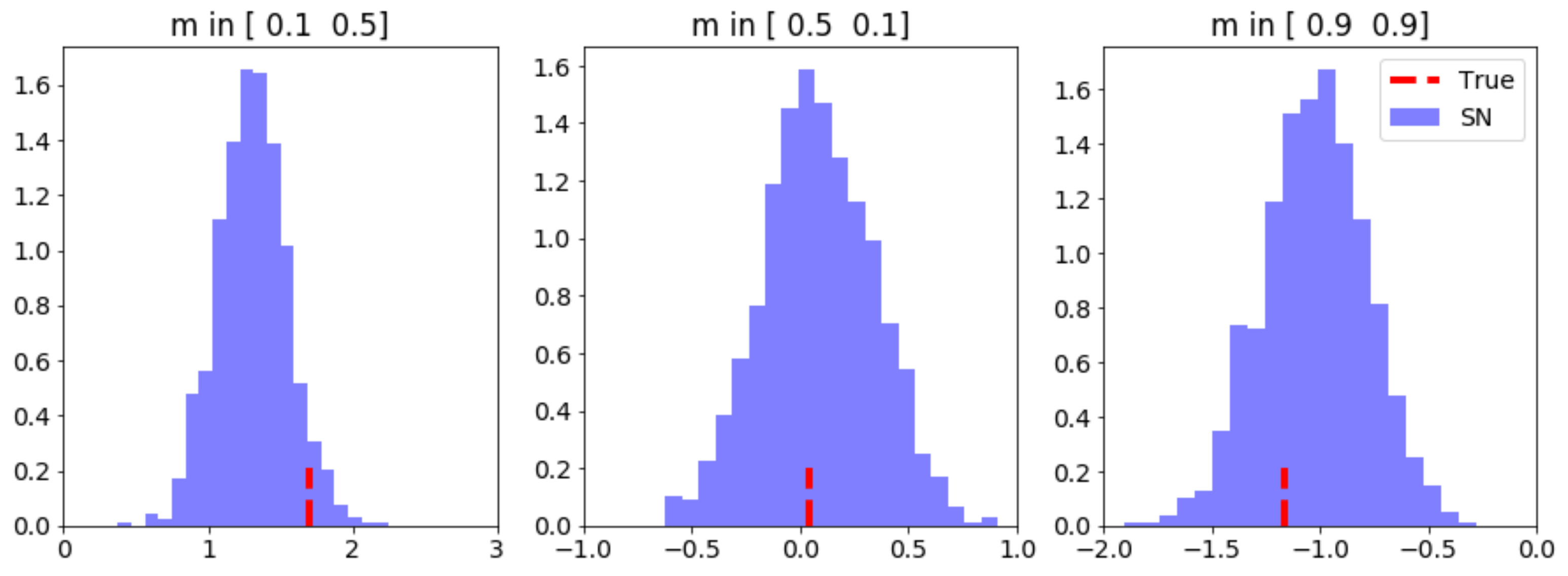}
    \includegraphics[width=0.95\textwidth]{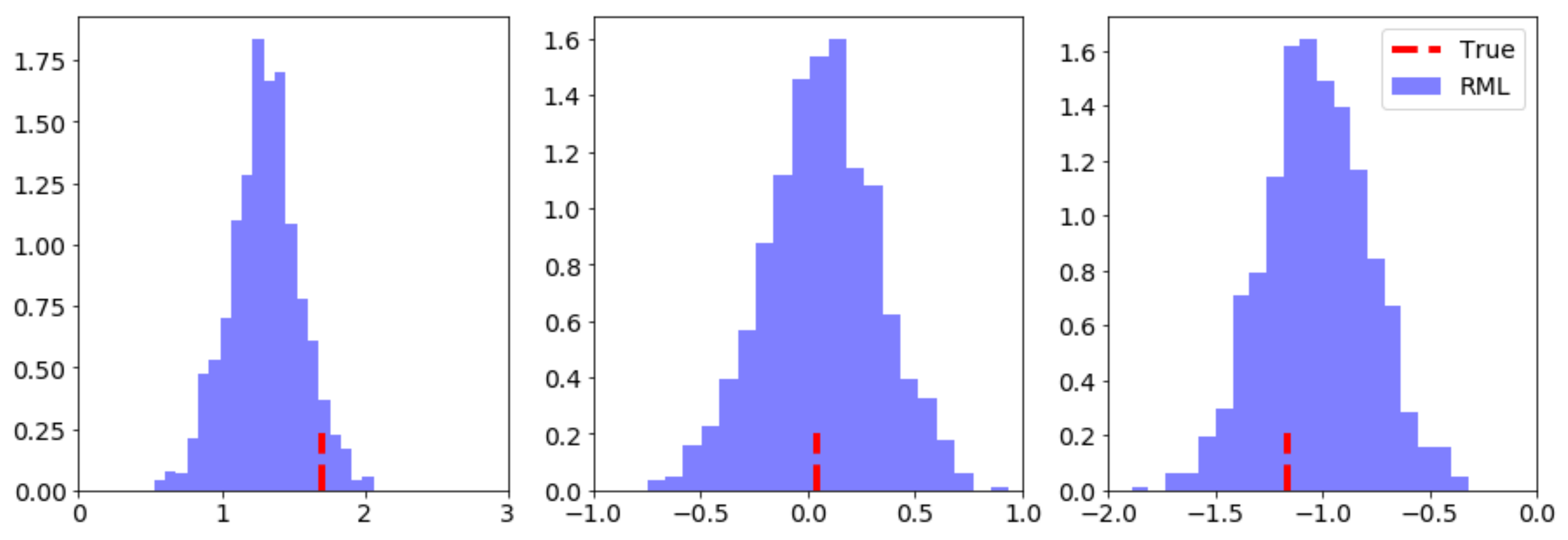}
    \includegraphics[width=0.95\textwidth]{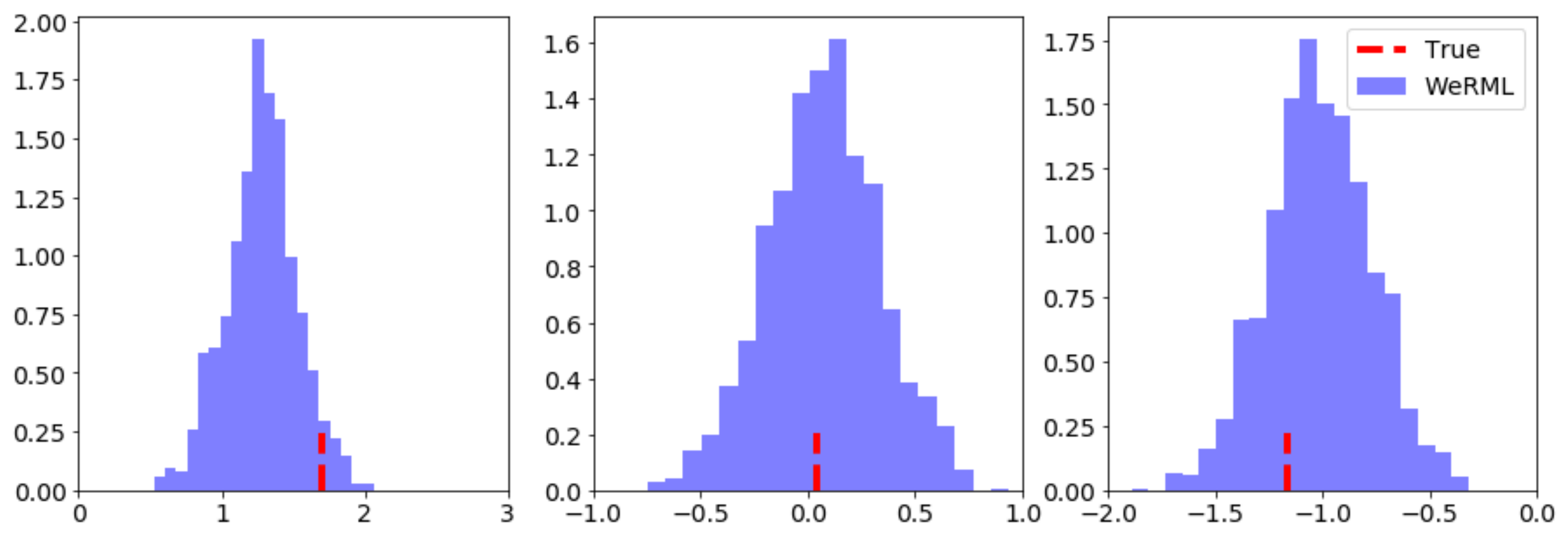}
  \caption{Density histograms of samples of the latent variable at locations (0.1,0.5), (0.5,0.1) and (0.9,0.9)  using SN (upper row), unweighted RML (middle row) and weighted RML (lower row).} \label{fig:Unw}
\end{figure}

In Fig.~\ref{fig:Unw}, we plot the distribution of sample values of the latent variable at three locations for which the true field has  values $m = 1.7$, $0.04$ and $-1.17$.
For this  example, the marginal distributions of samples  at observation locations from unweighted RML, SN and weighted RML  are all similar and approximately Gaussian.

\begin{figure}[tbp]
  \centering
    \includegraphics[width=0.95\textwidth]{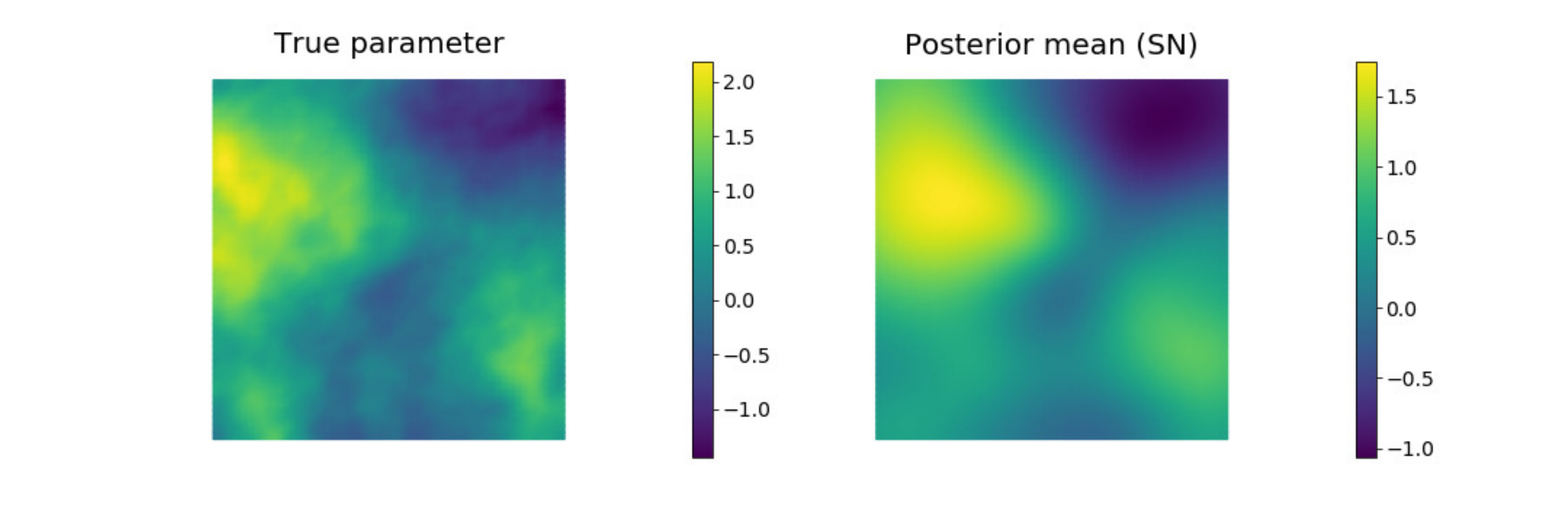}
    \includegraphics[width=0.95\textwidth]{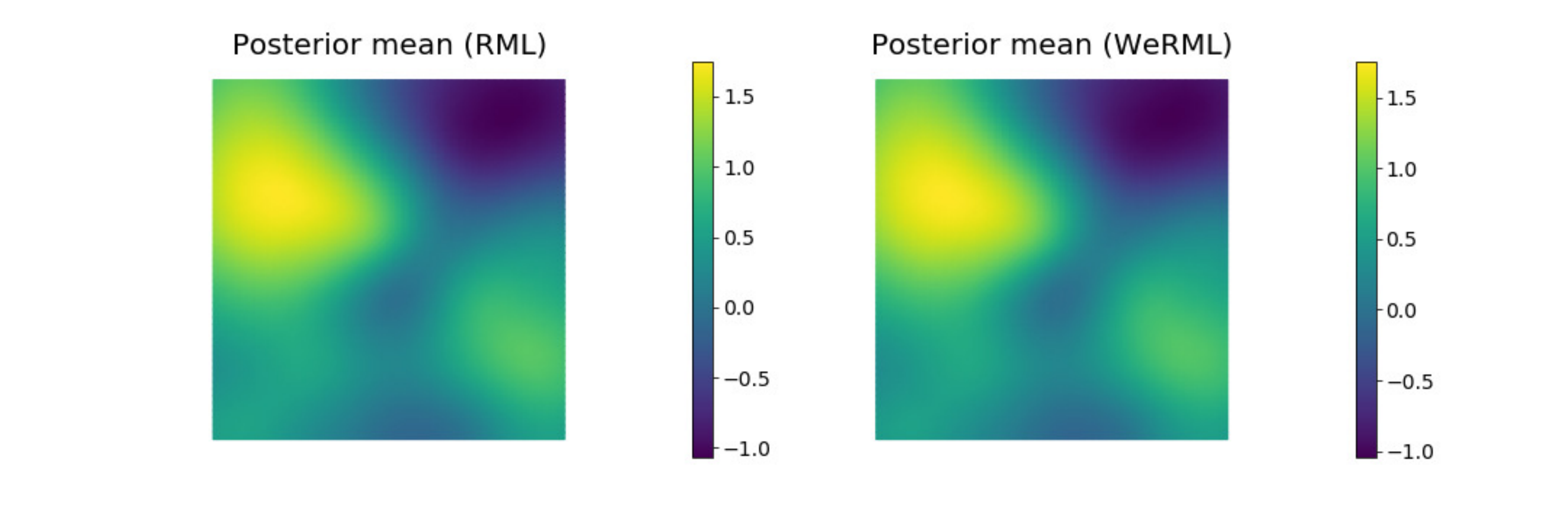}
  \caption{The true log-permeability field for Case 1 (upper left) and the posterior mean log-permeability fields from three approximate sampling methods: stochastic Newton (upper right), RML without weighting (lower left)  and RML with GN approximate weighting (lower right).} 
  \label{fig:loG}
\end{figure}

Estimates of the posterior mean of the log-perm field from three different sampling approaches are shown in Fig.~\ref{fig:loG}. As the permeability field is a monotonic function of the latent variable for this case, the estimated conditional means of the log-perm fields are similar  for the three methods.

\begin{figure}[tbp]
  \centering
\includegraphics[width=0.55\textwidth]{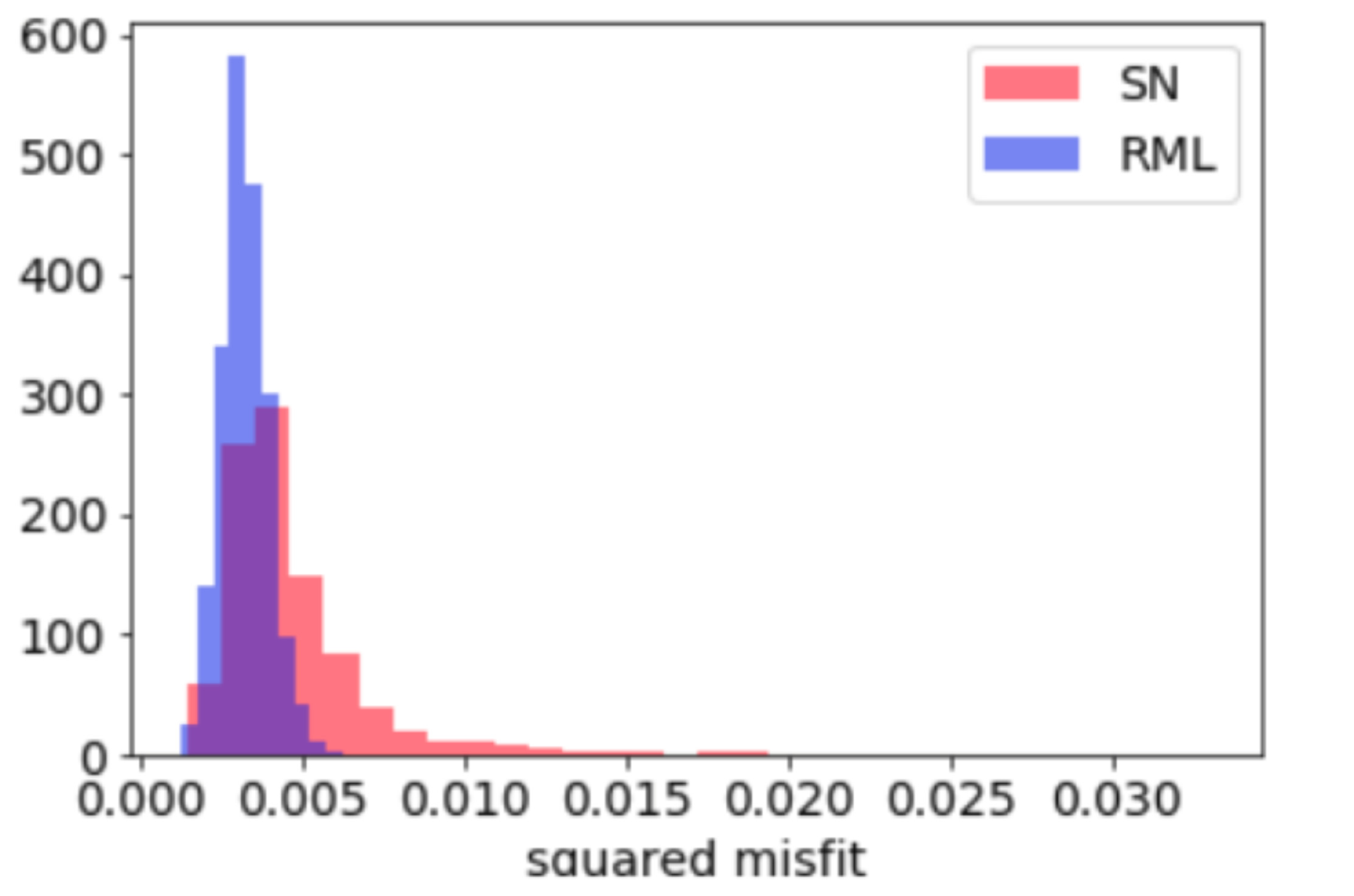}
  \caption{Distributions of squared data misfit from samples generated using SN and using weighted RML.}\label{fig:data_misfit_histograms}
\end{figure}

Despite the similarity of the mean fields  and the similarity of the estimates of the posteriori standard deviation from the three methods, the realizations from the three methods are not as similar in their ability to reproduce data (Fig.~\ref{fig:data_misfit_histograms}). The mean squared data misfit for weighted RML is $0.0032$, while the mean squared data misfit for stochastic Newton is $0.0048$. This is similar to the observation of Liu and Oliver \cite{liu:03c} who showed that the data mismatch of realizations generated from the posteriori mean and covariance in a 1D Darcy flow problem were much larger than data mismatch from MCMC or from RML.

\subsubsection{Case 2: log-permeability  is monotonic function of latent variable}

Although it is common to assume that permeability is log-normally distributed within a single rock type,  in many subsurface formations the distribution of permeability is largely controlled by `rock type'. In Case 2 we model the spatial distribution of rock types by applying a soft threshold to a latent Gaussian random field. With this transformation, values of $m < 1$ are assigned $\log \kappa \approx -2$ and values values of $m > 1$ are assigned $\log \kappa \approx 2$. One practical consequence of this transformation is that minimization of the objective function is more difficult. The second, more important consequence is that the nonlinearity in the neighborhood of the minimizers increases the variability in weights. So while in Case 1, approximately 90\% of the weights were between 0.0006 and 0.0015, in  Case 2 approximately 90\% of the weights fell between $10^{-4}$ and $9 \times 10^{-3}$ (Fig.~\ref{fig:Wemi_b}). The effective sample size for the 930 samples that converged successfully is also smaller in this case; $N\sbr{eff} = 179$ for an efficiency of about 19.2\%.

\begin{figure}[htbp!]
\centering
\subfloat[Case 2 (monotonic log-permeability)]{\label{fig:Case2_new_misfit}
\includegraphics[width=0.47\textwidth]{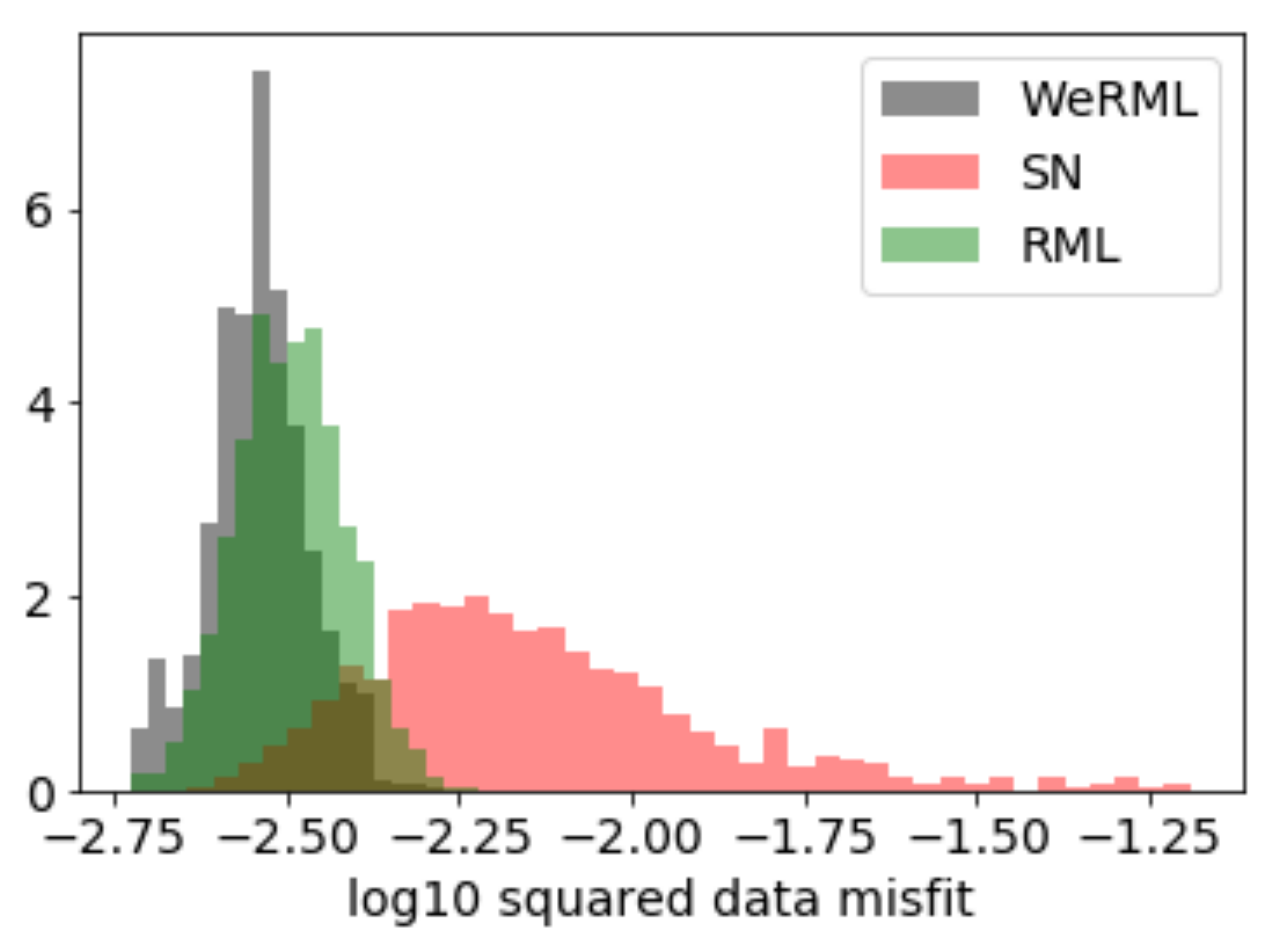}}
~
\subfloat[Case 3 (non-monotonic log-permeability)]{\label{fig:Case2_misfit}
 \includegraphics[width=0.45\textwidth]{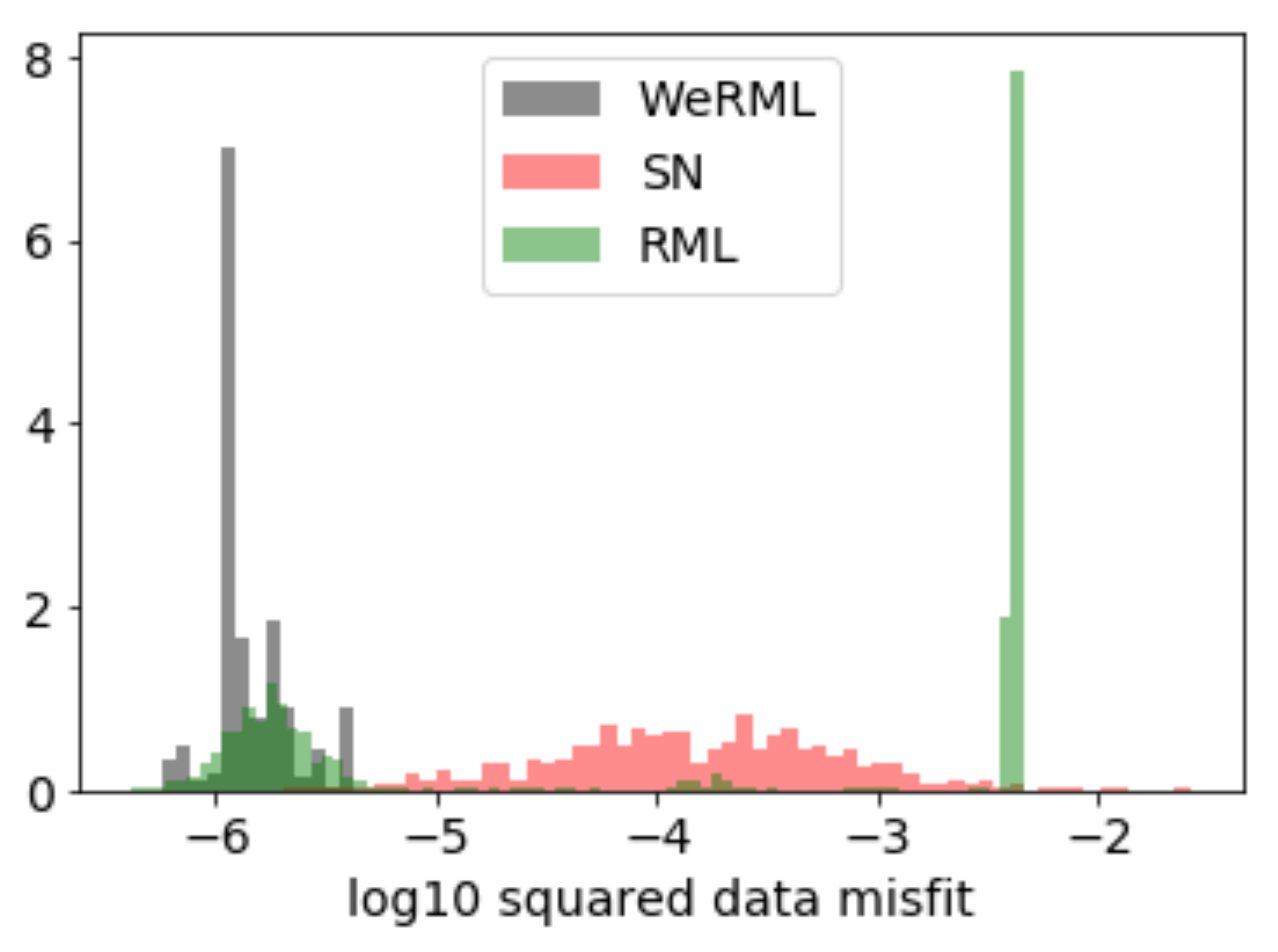}}
\caption{Compare distributions of squared data misfit for three sampling methods: RML, weighted RML (WeRML), and stochastic Newton (SN).}
\label{fig:compare_distributions}
\end{figure}

As the posterior distribution is not even approximately Gaussian, the mean and variance may not be the best attributes for judging the quality of the data assimilation. It is common in inverse problems to judge the quality of the realizations by the data misfit after calibration. The expected value of the mean squared data mismatch with observations is $0.0025$. The value computed from weighted RML ($0.0029$) is quite close to that value. In contrast, the value from SN realizations  ($0.0087$) is about $3.5$ times  larger than expected and the value from unweighted RML ($0.0033$) is slightly larger than weighted RML. The distributions of squared data misfits for the three methods are shown in Fig.~\ref{fig:Case2_new_misfit}.

\subsubsection{Case 3: log-permeability  is non-monotonic function of latent variable}\label{sec:non-monotonic}

\begin{figure}[tbp]
  \centering
    \includegraphics[width=0.95\textwidth]{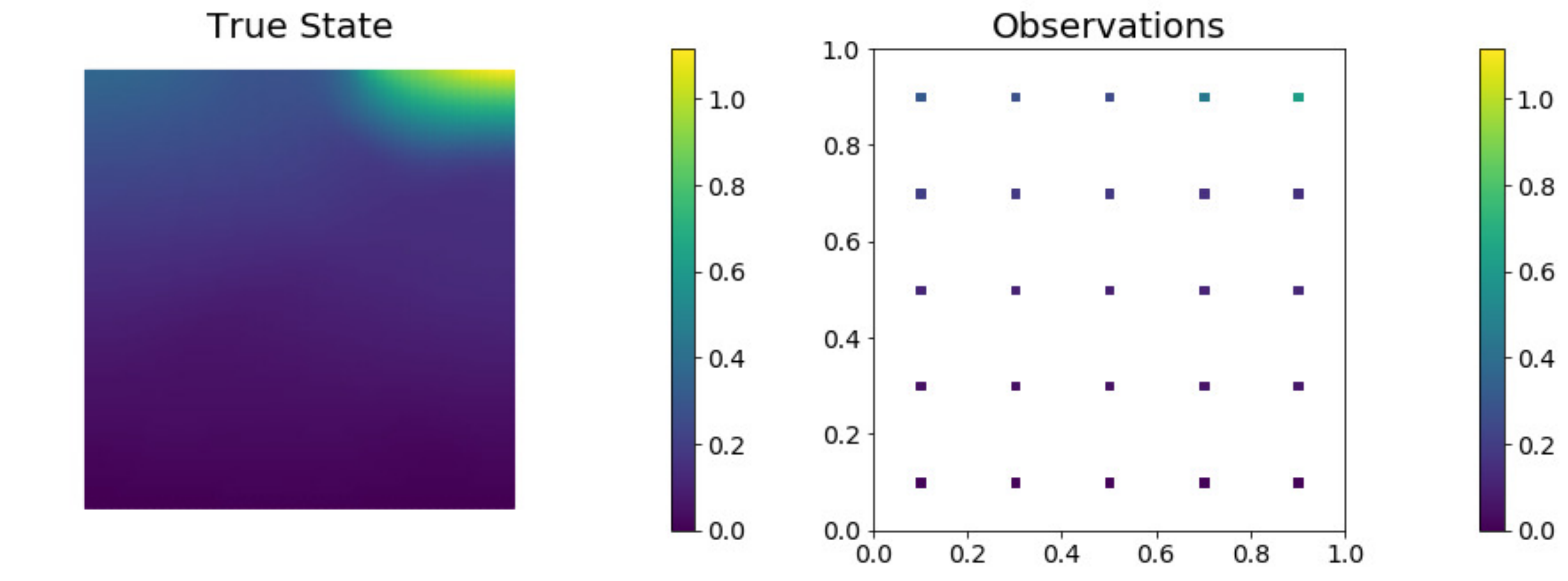}
  \caption{The true state (left) and observation locations (right)}\label{fig:Obs}
\end{figure}

Applying the transformation  \eqref{eq:threshold3},
we obtain a permeability field with a low permeability `background' and connected high perm `channels'. The `true' latent variable field and the corresponding true log-permeability field  are  shown in Fig.~\ref{fig:True}. The true pressure field, from which data are generated, and observation locations are plotted in Fig.~\ref{fig:Obs}.

In this example, the transformation from the Gaussian latent variable to the permeability variable is highly nonlinear and non-monotonic, so we should expect to encounter two problems: convergence to the minimizer will be slow \cite{liu:04a} and the algorithm is likely to converge to a local minimum that does not have large probability mass associated with it.

\begin{figure}[htbp]
  \centering
    \includegraphics[width=0.95\textwidth]{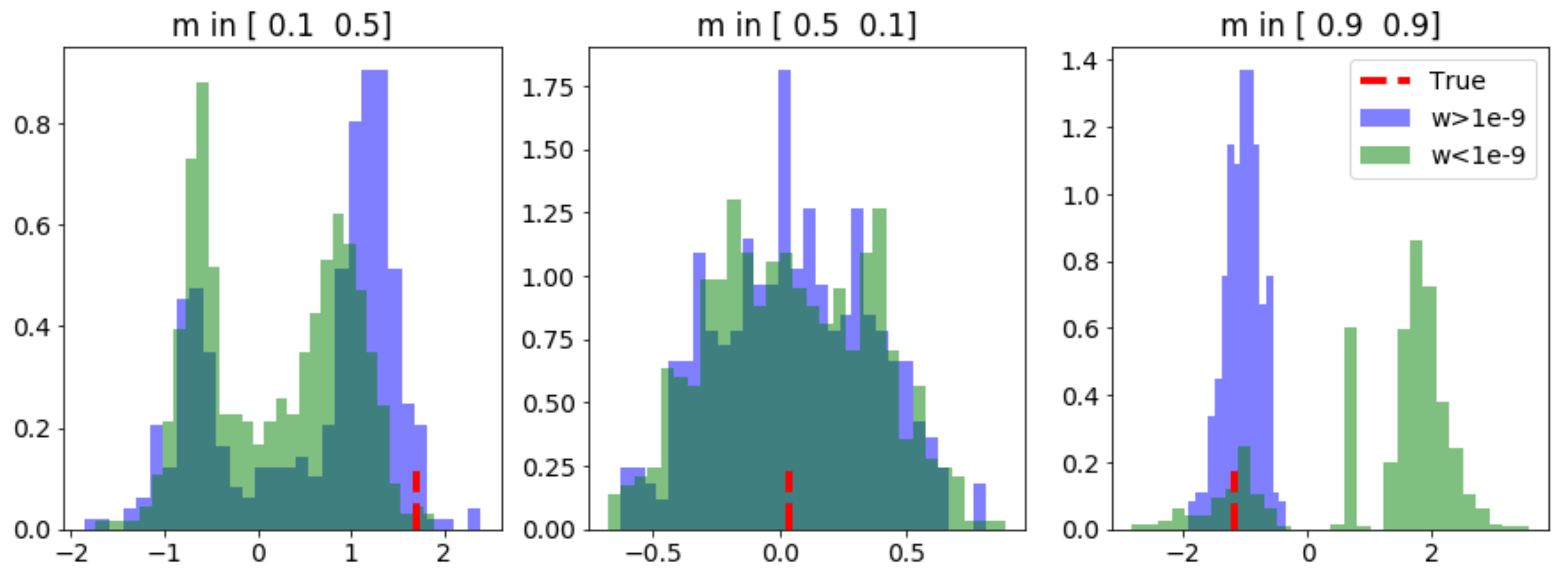}
    \includegraphics[width=0.95\textwidth]{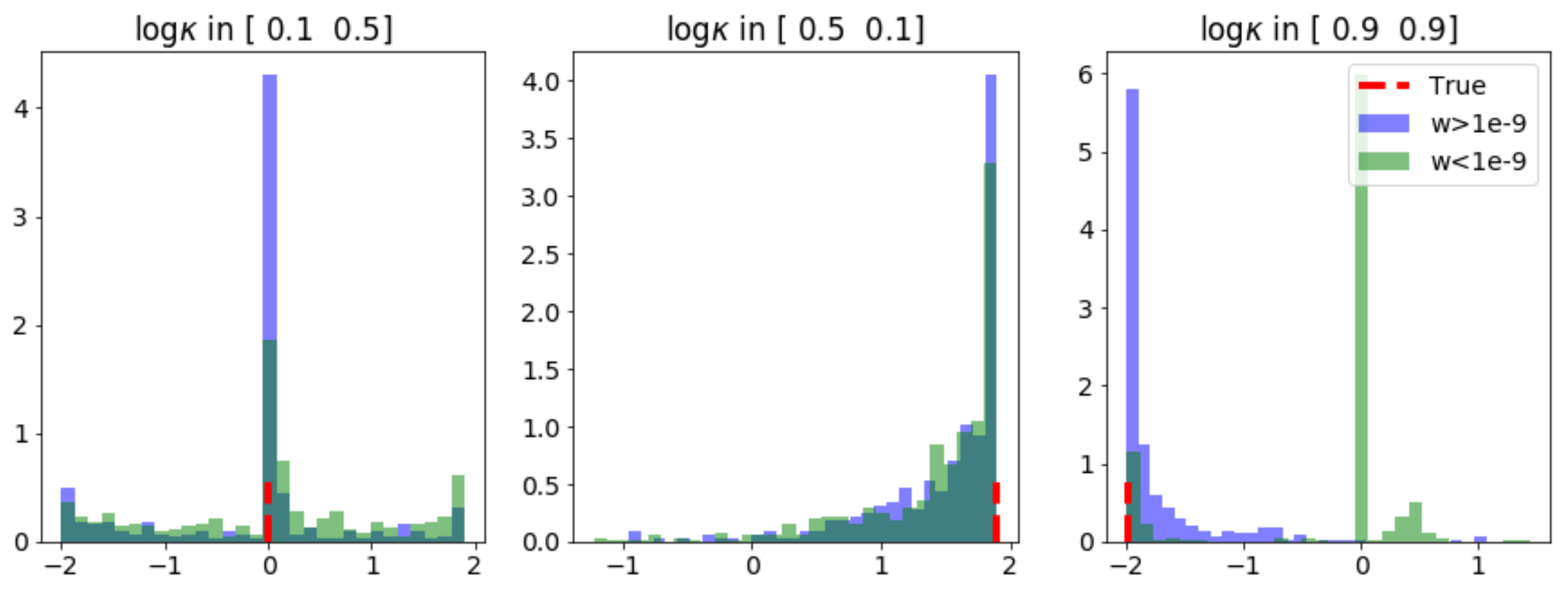}
 \caption{Density histograms of values of the minimizers of the objective functions at three observation locations (upper row) and corresponding values of the log-perm field (lower row). Colors separate minimizers into two groups  by weight.}\label{fig:Unw1}
\end{figure}

We focus on the distribution of samples that are obtained using a practice that could feasibly be applied to large-scale subsurface data assimilation problems if a gradient is available -- search for a single minimizer for each sample from the prior and use the Gauss-Newton approximation of the Jacobian to compute the weights. Here we have taken  1000 samples from the prior and performed 1000 corresponding minimizations, from which we obtained  $N_e=885$ samples with  successful termination and weights. Using
\eqref{eq:NEff}, we compute the effective sample size   $N_{\text{Eff}} \approx 14$.
Because of the dimensionality of the problem, it is not possible to completely characterize the posterior distribution for either the latent variables, or for the log-permeability. In order to gain some understanding, we examine the marginal distribution of the minimizers at three observation locations for which the true log-permeability values are approximately, -2, 0 and 2. The marginal distribution of unweighted minimizers (upper row Fig.~\ref{fig:Unw1}) is  bimodal at two of the observation locations.
The lower row of Fig.~\ref{fig:Unw1} (lower row) shows the corresponding distribution of unweighted log-permeability values at the same locations.
Although the spread of $m$ is fairly large at each of the observation locations, the spread of the log-permeability values are tightly centered on  $-2$, $0$ and $2$ in Fig.~\ref{fig:Unw1} (lower row). This is a consequence of the thresholding property of the  log-permeability transform.

The colors used for the density histograms in Fig.~\ref{fig:Unw1}  separate the samples into two groups: one in which $w>10^{-9}$ and the second for which $w<10^{-9}$. Note that in Fig.~\ref{fig:Unw1} (lower right), a substantial fraction of samples converged to a local minimizer with $\log \kappa \approx 0$ at $x,y =(0.9,0.9)$, which is far from the true value,  $\log \kappa\spr{true}\approx -2$. Many of the clearly erroneous minimizers are easily eliminated, however, by the low weights.  Note that the inefficiency of the sampling in this case  is a result  of the non-monotonic nature of the permeability transform -- to get from  $\log \kappa =0$, which is a local minimizer, to the correct root  $\log \kappa =-2$ it is necessary to pass through the point  $\log \kappa =2$. In the histograms, the samples with high weights (blue bars) are generally close to the true values (red dashed line), which is as we expected.

\begin{figure}[tbp]
  \centering
    \includegraphics[width=0.95\textwidth]{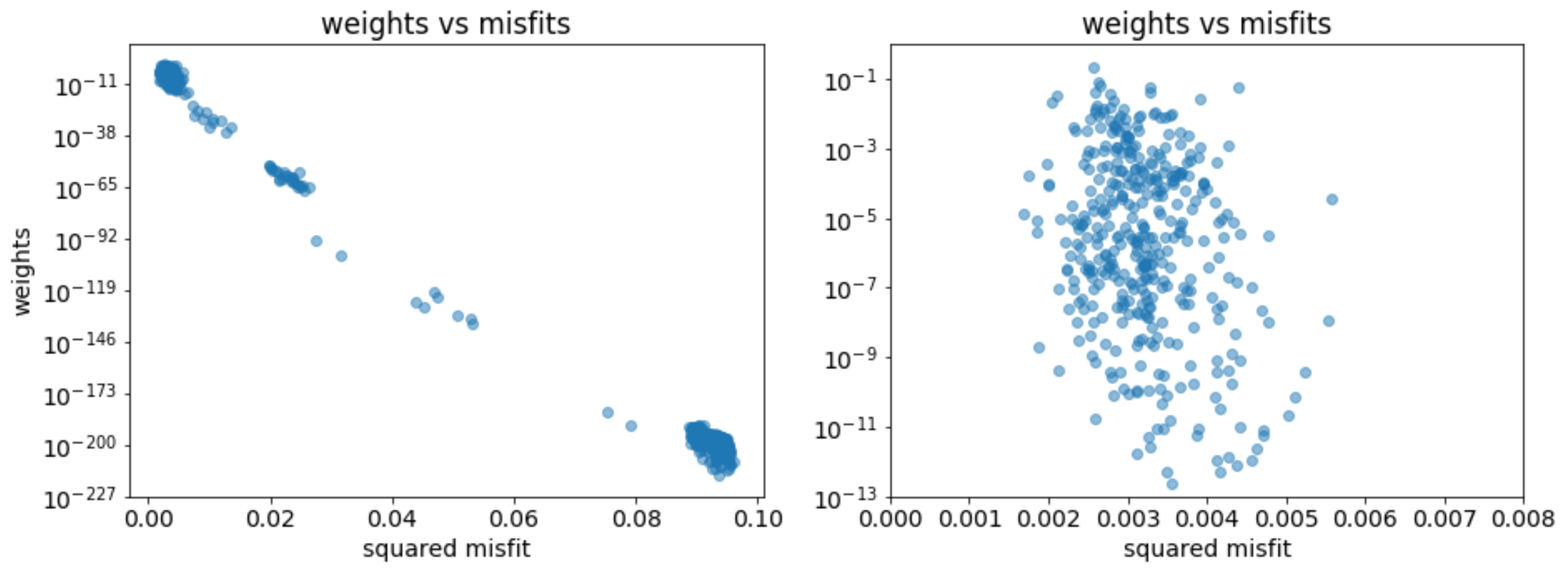}
  \caption{The weights vs misfits for the minimizers of the stochastic cost function.}\label{fig:Wemis1}
\end{figure}

As in the previous case, the expected value of the squared data misfit with the actual observation is approximately $0.0025$ at the global minimizer of the stochastic cost function.
This is close to the values that are obtained in the best minimizations (Fig.~\ref{fig:Wemis1} (right)). In this example, however, many of the minimizations converged to minimizers with much larger misfit values (Fig.~\ref{fig:Wemis1} (left)). The minimizations with large misfit values result from convergence to local minima in the cost function. The number of local minima with large numbers of samples appears to be relatively small as a result of the large correlation range for the latent permeability variable. For this example, it appears that many of the unwanted local minimizers  could be eliminated either through the weighting, or through the magnitude of the squared data misfit.

\begin{figure}[tbp]
  \centering
  \includegraphics[width=0.95\textwidth]{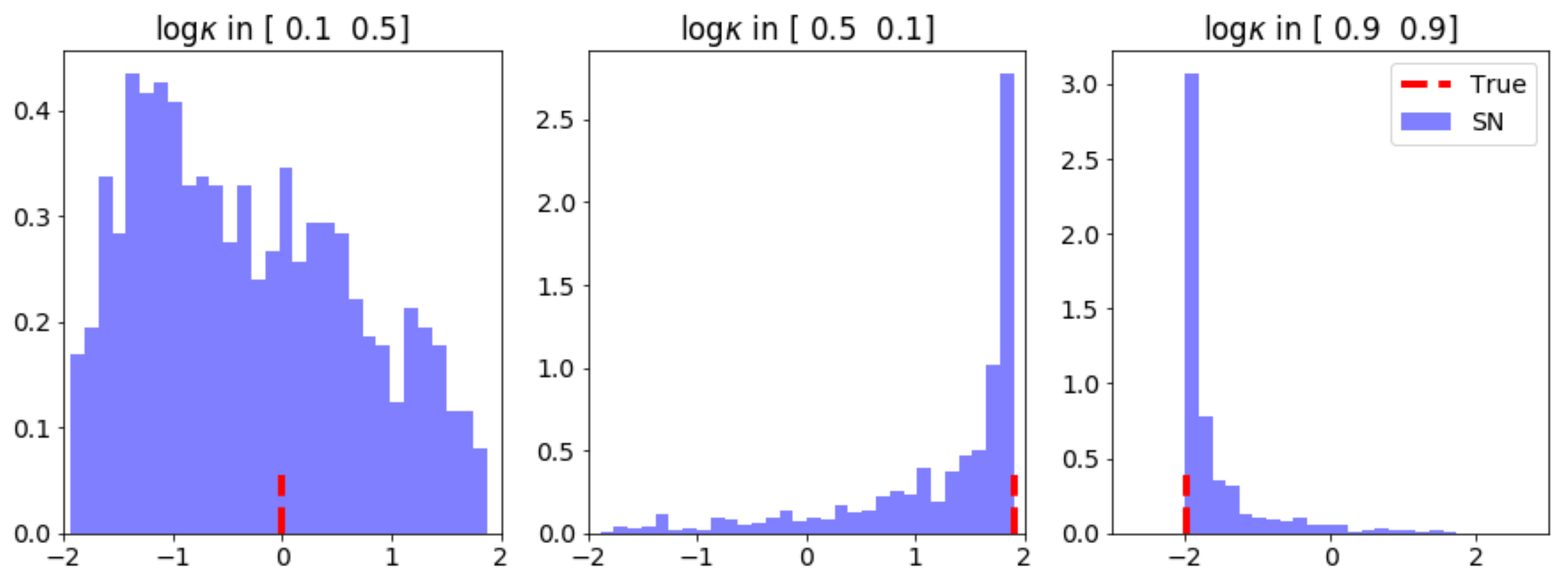}
    \includegraphics[width=0.95\textwidth]{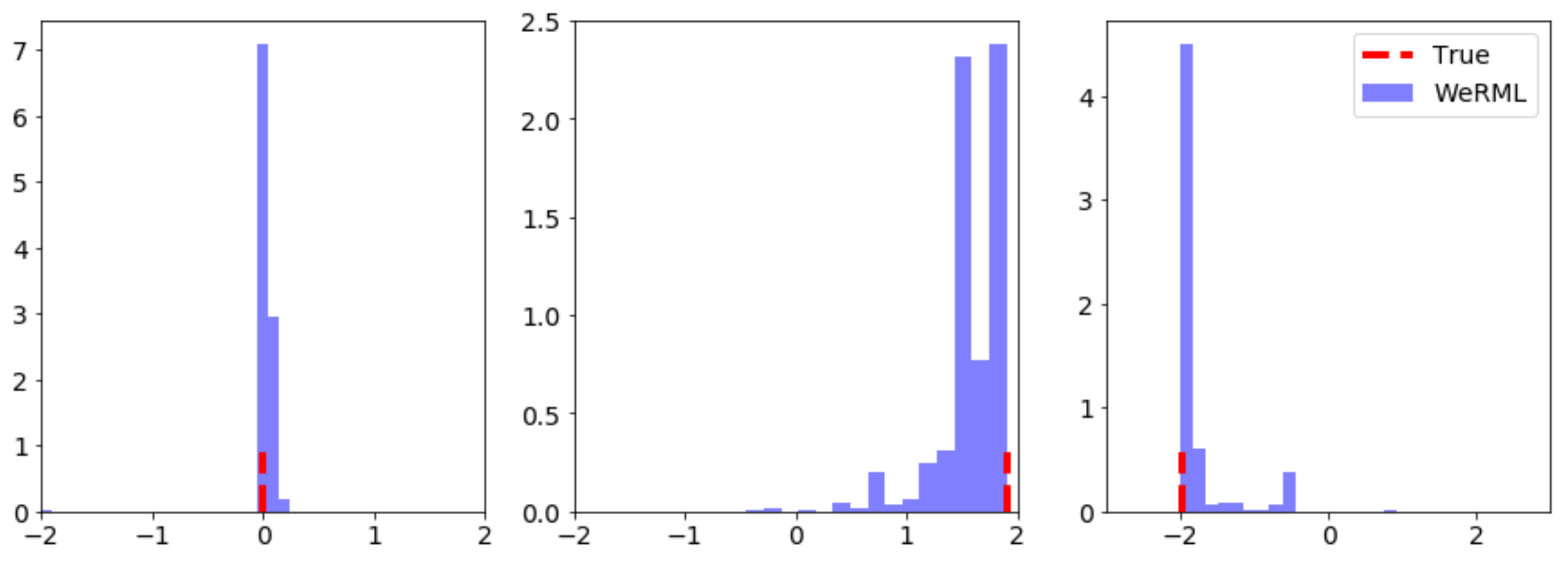}
  \caption{The histograms of values of the log-perm field at three observation locations by using the SN (upper row) and weighted RML (lower row) methods.}
  \label{fig:SN-RML}
\end{figure}

Figure~\ref{fig:SN-RML} compares distributions of the values of the log-permeability field obtained by the SN and the  weighted RML methods  at three pressure observation locations.  For consistency, the same sample size was used for both methods.
The true values of the log-permeability field at the observation locations are close to $-2$, $0$ and $2$ (shown as small red dots).
When using the weighted RML samples to approximate the distribution, the samples  are seen to be concentrated close to the true values at the three observation locations (Fig.~\ref{fig:SN-RML} (lower)).
Because of the highly nonlinear transformation of the permeability field, the SN method was unable to provide a good approximation to the true posterior distribution, although it did provide plausible estimates of uncertainty in  $\log \kappa$ at two of the locations. At the third location (upper left  in Fig.~\ref{fig:SN-RML}) the uncertainty in log-permeability was completely misrepresented.

\begin{figure}[tbp]
\centering
    \includegraphics[width=0.95\textwidth]{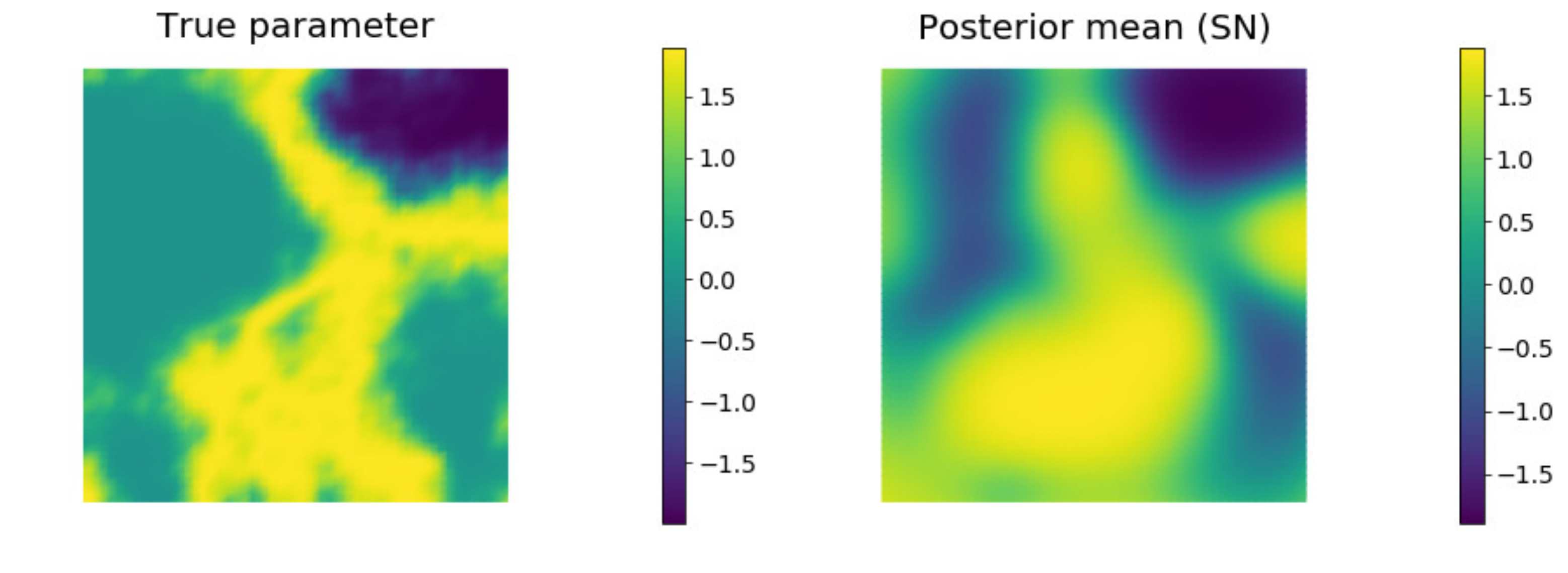}
    \includegraphics[width=0.95\textwidth]{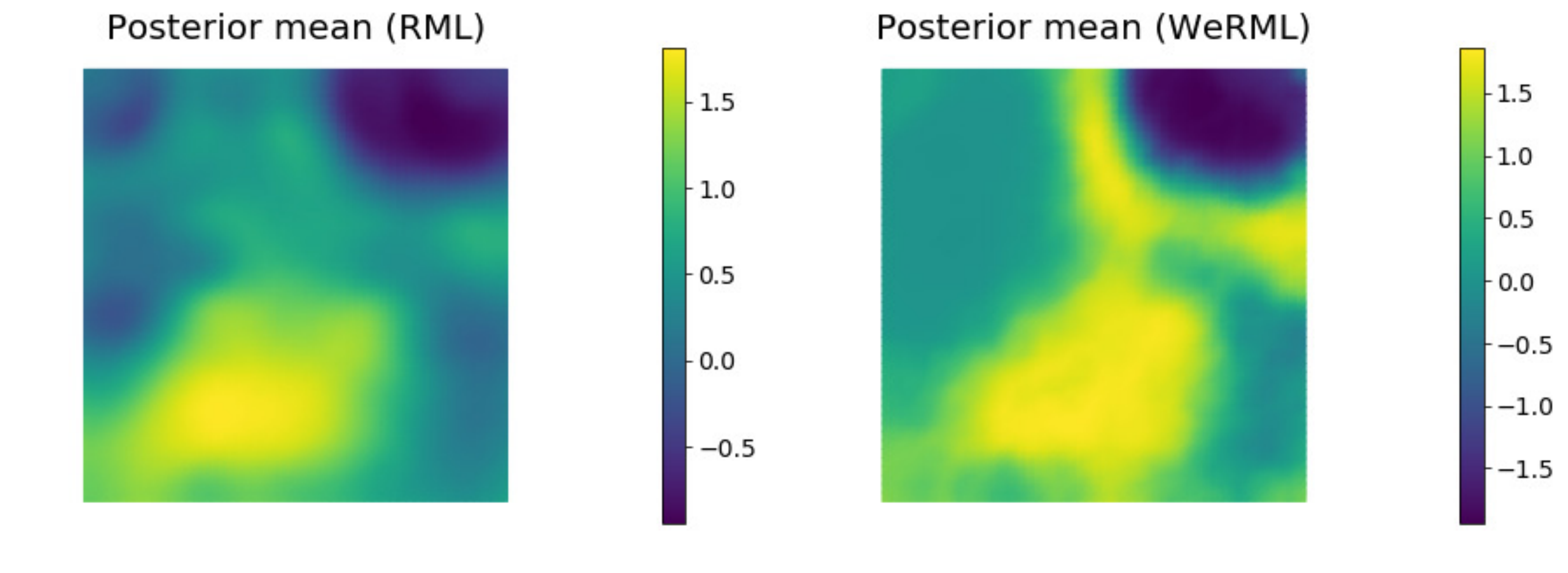}
  \caption{ The true log-permeability field for Case 3 (upper left) and the posterior mean log-permeability fields from three approximate sampling methods: stochastic Newton (upper right), RML without weighting (lower left)  and RML with GN approximate weighting (lower right).}   \label{fig:loG1}
\end{figure}

Although the marginal distribution of the log-permeability field is multimodal at some locations before and after data assimilation, the mean of log-permeability is useful for qualitatively gauging the quality of the data assimilation. In Fig.~\ref{fig:loG1}, we compare the MAP point of SN, and the unweighted and weighted posterior means of RML with the true log-permeability field. All three methods provide reasonable characterization of the mean log-permeability in the upper right area of the grid, but only the weighted RML method adequately characterizes log-permeability on the left side. It appears that it would be risky to use results from RML without weighting in this case.

\begin{figure}[htp]
\centering
     \includegraphics[width=0.95\textwidth]{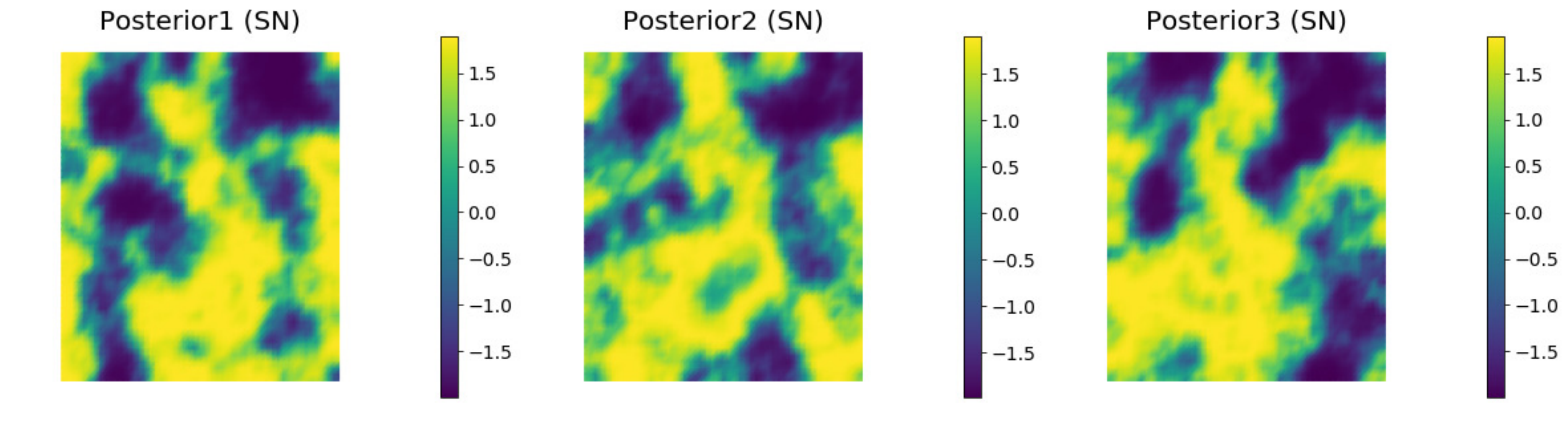}
    \includegraphics[width=0.95\textwidth]{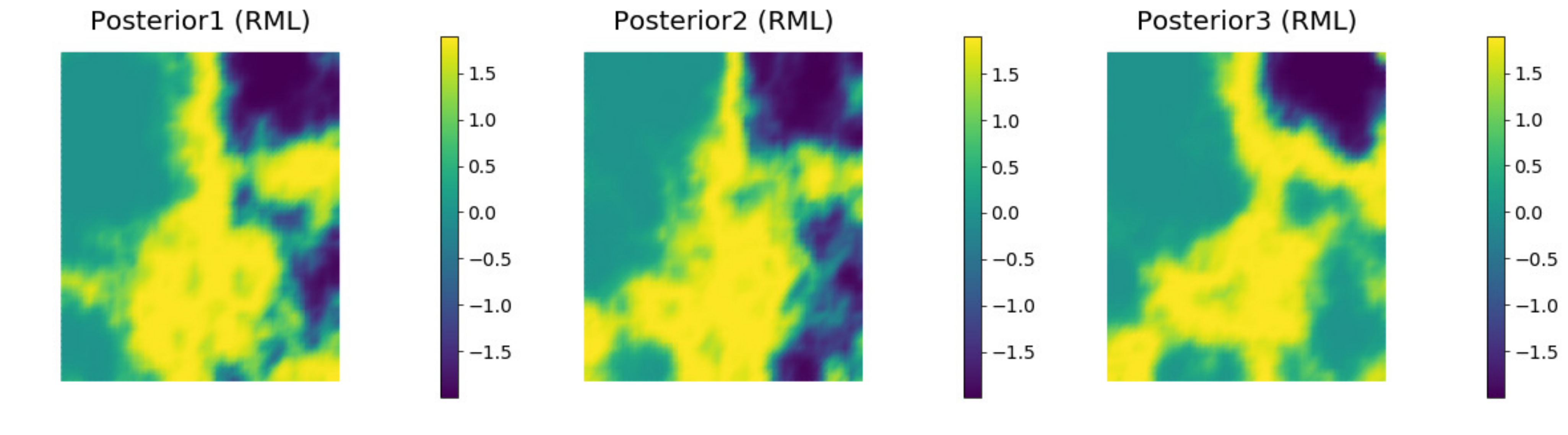}
  \caption{Three posterior log-permeability samples from the stochastic Newton method (top row) and the three weighted RML samples with largest weights (bottom row). Samples can be compared with the true log-permeability distribution in Fig.~\ref{fig:loG1}.}\label{fig:Pos}
\end{figure}

As the posterior distribution is multi-modal, the mean and variance may not be the best attributes for judging the quality of the data assimilation. Two qualitative criteria are often used in practice. First, it is common to judge the quality of the realizations by the data misfit after calibration. The expected value of the mean squared data mismatch with observations is $0.0025$. The value computed from weighted RML is quite close to that value, $0.0029$. In contrast, the value from SN realizations is about $10$ times larger than expected ($0.0264$) and the value from unweighted RML is larger still ($0.0486$). The distributions of squared data misfits for the three methods are shown in Fig.~\ref{fig:Case2_misfit}.  Second, the samples themselves can be examined qualitatively for `plausibility' -- do they look like samples from the prior?
Fig.~\ref{fig:Pos} shows RML samples with largest weights (bottom row) and the corresponding posterior samples for SN (top row).
In this case, the weighted and unweighted RML samples look plausible, but the samples from SN do not.

\subsection{Discussion of porous flow results}

The porous flow examples were chosen to be large enough that a na{\"i}ve approach to particle filtering in which particles are sampled from the prior and then weighted by the likelihood would suffer from the curse of dimensionality and all the weight would fall to a single particle. The dimension of the model space ($2601$ discrete parameters) was also large enough that computation of the Jacobian of the transformation  from the prior distribution to the distribution of critical points would be challenging.

Exact sampling of the posterior distribution using this methodology requires either computation of all critical points of the objective function, or random sampling of all critical points. In the porous media flow examples, however,  we were unable to locate any maximizers for the objective function even for Case 3 in which the objective function had many local minima. As a consequence, it appears that searching only for minimizers is a robust approximation in high dimensions. Because we searched only for minimizers, it also appears that the Gauss-Newton approximation  of the Jacobian gave useful approximations. The terms that must be computed are then very similar to terms that are computed in Gauss-Newton minimization of the cost function. It was possible to compute an inexpensive estimate of the  determinant of the Gauss-Newton approximation of the Jacobian using eigenvalues of the Hessian at the minimizers. Because the dimensions of the data space was relatively small, it was also possible to estimate the determinant of the Jacobian using the adjoint system. For Case 1, the estimates from the two approaches were similar, but the differences increased as the nonlinearity increased. 
To evaluate the quality of the sampling from the various methods,  we used weights computed using the adjoint system. 
The weights were more variable when low rank approximations were used. In that case, it was useful to account for the model error by inflating the value of $C_d$ used for computing weights.

The degree of nonlinearity in the transformation from parameter to log-permeability had a strong effect on the effective efficiency of the minimization approach for sampling. The least nonlinear example (Case 1) had an effective sample efficiency of 84\% while the most nonlinear example (Case 3) had an effective sample efficiency of 1.6\%.  It appears that the low efficiency in Case 3 was largely a result of the prevalence of many local minima in the objective function, many of which were characterized by large data mismatch and very small weights.

When the posterior distribution had a single mode as in Case 1, the distribution of residual errors in the data mismatch was quite small and the correlation between weight and data mismatch was correspondingly small ($r=-0.086$). In that case, the data mismatch would not have provided a useful proxy for weighting. In Case 2, the nonlinearity was greater but it appears that the posterior distribution was still uni-modal. The weights did correlate with data mismatch in that case ($r=-0.485$) but it appears that the skewness of the distribution may have been the largest reason for the decrease in effective sample efficiency.  Finally, in Case 3, the transformation from the model parameter to log-permeability was non-monotonic and the posterior distribution was characterized by a large number of local minima. Here, the correlation between  importance  weight and  data mismatch was almost perfect ($r=-0.999$) and the data mismatch could serve as a useful tool for eliminating samples with small weights.

Although we did not compare the distribution of samples from weighted RML with methods such as MCMC, we did compare with the stochastic Newton method because it is a practical and scalable method for approximate sampling in high dimensions. For Case 1, which appears to be unimodal, the mean log-permeability fields from SN and RML (both weighted and unweighted) were visually similar. For Case 3, the  mean permeability fields from SN and weighted RML are less similar, the data mismatches from SN are substantially larger, and the samples are visually less plausible.

The documented cost of the three considered methods are substantially different. 
The computational complexity of the stochastic Newton method stems from a single minimization to compute the MAP and the cost to generate samples from a low-rank  Gaussian approximation of the posterior. While the sampling step contributes to the cost of the stochastic Newton, it is dominated by minimization of the objective function. For the Darcy flow example, the cost to generate 1000 approximate samples varied from 13 seconds for the log-normal case to 56 seconds for the non-monotonic case.\footnote{Timing should be considered illustrative, but for reference all results were obtained on a computer with a $\text{i7-5500U@2.40GHz}\times 4$ processor with 7.5  GiB memory and a 64-bit operating system. } Note that this increase is a result of the varying number of iterations required for the minimizer to converge for the different settings. 
In case of the RML method, the cost to generate $N_e$ realizations is dominated by the cost to perform $N_e$ minimizations with different cost functions. Henceforth the computational complexity for RML can be expected to be approximately $N_e$ times greater than the cost for stochastic Newton method. Indeed, in our examples, the run time required to generate 1000 samples from RML was approximately 1000 times greater, varying from 13000 seconds for the log-normal case to 47000 seconds for the non-monotonic case. 
For weighted RML, there is an additional cost incurred in the computation of the weights. Although several of the terms in the weights can be obtained at low cost through the same low-rank approximations that were used for the Hessian, we chose tosolve the adjoint system $N_d$ times to compute the Jacobian of the data for computation of $V\inv$. The additional cost for computing the weighting is thus dominated by the cost of running the simulator an additional  $N_d$ times for each realization. Because the adjoint system was solved to compute weights, the  cost of computing weights varied from 24000 seconds to 35000 seconds for 1000 samples, which was similar to the cost of the minimization.
All computational costs, including the cost of minimization, could be reduced through careful modification of the algorithms. In particular, the efficiency of the weighted RML could be improved by tempering the objective function at early iterations to avoid  convergence to local minima with small weights. Also, the cost of computation of the weights  could be reduced by using a low-rank approximation of $V$ as in the ensemble Kalman filter.

\section{Summary}

We have presented a method for sampling from the posterior distribution for inverse problems in which the prior distribution of model variables and measurements errors are Gaussian. Although the method is highly efficient when the posterior distribution is also approximately Gaussian, the target application is to problems in which the posterior distribution is multimodal  -- situations in which Gaussian approximations of the posterior distribution are inappropriate.
Because of the requirement that the prior distribution be Gaussian, this method then is probably more appropriate for parameter estimation problems  than for state estimation problems in which the prior may be multimodal as a result of nonlinear dynamics \cite{morzfeld:19}.
The method is similar to the method of randomized maximum likelihood or randomized maximum a posteriori in that samples are generated from the prior distribution and then moved  to regions of high probability. Instead of solving only for minimizers of a stochastic cost function, however, the method samples correctly when \emph{all} critical points are sampled and weighted, or when the critical points are \emph{randomly sampled} and weighted.
This procedure sampled correctly in small multimodal and skewed toy problems for which it was possible to compute all critical points. In those cases, it was also possible to obtain good approximate sampling using only minimizers of the cost function and a relatively inexpensive Gauss-Newton approximation of the particle weights.

The toy problems showed that the weights on maximizers are generally small and cannot be computed accurately using the Gauss-Newton approximation. Consequently a practical approach for larger inverse problems is to compute only the minimizers of the cost function and use approximate weights. We  showed that the weights can be computed from low-rank approximations of the Hessian evaluated at the minimizers. This approach was applied to three porous media flow examples. In the first example the posterior pdf appears to be unimodal. In that case, the spread in the weights was small and the effective sample efficiency was 84\%. The distribution was not visibly different from samples obtained using a less expensive low rank approximation of the Hessian, but the quality of the match to the data was significantly better.

The flow example with a non-monotonic transformation to log-permeability was much more difficult. In this case,  the cost function was characterized by a large number of local minimizers with small probability mass. Many of the local minimizers could be easily rejected on the basis of either low weights or unexpectedly large mismatch with observations. Because of the difficulty of converging to minimizers with large weights, the variance of the weights was large and the effective sampling efficiency in this case was low, approximately 1.6\%.
We emphasize, however, that this example was chosen to be extremely nonlinear to test the ability to sample the multimodal posterior distribution for a moderately large model with thousands of parameters. The weighted mean of the samples from RML in this case, provided a good approximation to the true permeability distribution and the weighted data mismatches were close to the expected value.

\section*{Conflicts of interest}

The authors declare that there are no conflicts of interest.

\section*{Availability of data}

No data are used in the manuscript.

\section*{Code availability}

Selected Python codes used in the preparation of this manuscript are available at \url{https://bitbucket.org/JanadeWiljes/workspace/projects/RBPS.}

\section*{Funding}
For this work, Yuming Ba was supported by the China Scholarship Council. Dean Oliver was supported by the  NORCE Norwegian Research Centre cooperative research project ``Assimilating 4D Seismic Data: Big Data Into Big Models'' which is funded by industry partners Aker BP, Equinor, Lundin Norway, Repsol, and Total, as well as the Research Council of Norway through the Petromaks2 program. Jana de Wiljes and Sebastian Reich have been partially funded by Deutsche Forschungsgemeinschaft (DFG) - Project-ID 318763901 - SFB1294.

\appendix

\section{Nomenclature}

\begin{table*}[ht]
\begin{center}
\begin{tabular}{ll}
$d^{\rm o}$ & observations  \\
$g$ & forward map  \\
 $m$ and $\delta$ & samples of target distribution  \\
$m_i$ & $i$th sample of posterior distribution  \\
$m^\ast$ & unknown reference parameter \\
$m'$ and $\delta'$ & samples of Gaussian $q_{M'\Delta'} (m',\delta')$   \\
$m_\text{MAP}$ & optimizer of cost functional  \\
$\bar m$ & mean of Gaussian prior  \\
$n(m')$ & total number of critical points  \\
$n(z')$ & cardinality of $\mathcal{M}_{z'}$  \\
$p_{M\Delta}(m,\delta)$ & target density \\
$q_{M'\Delta'} (m',\delta')$ & proposal distribution  \\
$u(x)$ & pressure   \\
$w$ & weights  defined in \eqref{eq:weight-RML} \\
$z$ & tuple of samples $m$ and $\delta$   \\
$z'$& tuple of samples $m'$ and $\delta'$   \\ \hline
$A_0$, $A_1$, $A_2$ & normalization constants \\
$C_M$ & covariance matrix of Gaussian prior  \\
$C_D$ & covariance matrix of Gaussian observation error  \\
G & differential operator of $g$  \\
$H_\text{map}$,  $H_\text{misfit}$ & Hessian matrices   \\
$J$ & Jacobian determinant  \\
$L(m)$&log likelihood function \\
$\mathcal{M}_{z'}$ &set mapping to $z'$ via $\Psi$  \\
${\rm N}(\cdot,\cdot)$ & Gaussian distribution \\
$N_d$ & dimension of observational space  \\
$N_e$ & number of samples   \\
$N_m$ & number of model parameters  \\
$N_s$ & number of samples of prior  \\
$N_{\text{Eff}}$ &effective sample size  \\
$Q$ & square root of $C_M$  \\
$V$ & auxiliary variable defined in \eqref{eq:V}  \\ 
\hline
$\epsilon$ & observation error  \\
$\eta(m)$ & auxiliary variable defined in \eqref{eq:eta}  \\ 
$\kappa\spr{true}$ & reference permeability field  \\
$\kappa (x)$ & permeability field  \\
$\lambda_i$ & eigenvalues \\
$\pi_D(d^{\rm o})$ & normalization constant  \\
$\pi_M(m|d^{\rm o})$ & posterior distribution  \\
$\sigma_d$& standard deviation of measurement variables  \\
$\sigma_m $ &standard deviation of parameter variables  \\ \hline
$\Lambda$ & diagonal matrix of eigenvalues  \\
$\Psi$ & map between samples  \\
\end{tabular}
\end{center}
\caption{Notation used throughout the manuscript.}
\label{tab:Notations}
\end{table*}

\end{document}